\documentclass{gtart_h}  

\def\ifplaintex{\expandafter\ifx\csname documentclass\endcsname\relax}

\ifplaintex 
\else
\fi

\expandafter\ifx\csname beginpicture\endcsname\relax
\expandafter\ifx\csname documentclass\endcsname\relax
\input pictex \else
\input prepictex \input pictex \input postpictex \fi\fi

\def\gt{{\mathsurround=0pt\it $\cal G\mskip-2mu$eometry \&\ 
$\cal T\!\!$opology}}        

\def\gtp{{\mathsurround=0pt\it $\cal G\mskip-2mu$eometry \&\ 
$\cal T\!\!$opology $\cal P\!$ublications}}  

\def\lognumber#1{\def\thelognumber{#1}}
\def\volumenumber#1{\def\thevolumenumber{#1}}
\def\papernumber#1{\def\thepapernumber{#1}}
\def\volumeyear#1{\def\thevolumeyear{#1}}

\def\pagenumbers#1#2{\def\startpage{#1}\def\finishpage{#2}}
\def\published#1{\def\publishdate{#1}}
\def\proposed#1{\def\theproposer{#1}}
\def\seconded#1{\def\theseconders{#1}}
\def\received#1{\def\receiveddate{#1}}

\def\accepted#1{\def\accepteddate{#1}}

\def\asciiaddress#1{\def\theasciiaddress{#1}}
\def\asciiemail#1{\def\theasciiemail{#1}}

\long\def\asciiabstract#1{\long\def\theasciiabstract{#1}}

\let\\\par\let\thelognumber\relax
\let\thevolumenumber\relax\let\thepapernumber\relax
\let\thevolumeyear\relax\let\thesamplenumber\relax\let\startpage\relax
\let\finishpage\relax\let\publishdate\relax\let\receiveddate\relax
\let\reviseddate\relax\let\accepteddate\relax\let\theasciititle\relax
\let\theasciiauthors\relax\let\theasciiaddress\relax
\let\theasciiabstract\relax
\let\theasciiemail\relax\let\theshortauthors\relax\let\theshorttitle\relax

\long\def\maketitlep{   

\count0=\startpage

\gt\hfill      
\beginpicture
\setcoordinatesystem units <0.33truein, 0.33truein> point at 2.2 0.9
\setplotsymbol ({$\cal G$})
\plotsymbolspacing=9truept
\circulararc 315 degrees from 0 1 center at 0 0
\setplotsymbol ({$\cal T$})
\circulararc 315 degrees from 1 -1 center at 1 0
\endpicture
\break
{\small\ifx\thesamplenumber\relax 
Volume \else Sample
\fi\thevolumenumber\ (\thevolumeyear)
\startpage--\finishpage\nl
Published: \publishdate}
\vglue 0.5truein plus 0.4fil minus 0.1truein

{\parskip=0pt\leftskip 0pt plus 1fil\def\\{\par\smallskip}{\ifplaintex\large
\else\Large\fi\bf\thetitle}\par\medskip}   

\vglue 0pt plus 0.1fil 

{\parskip=0pt\leftskip 0pt plus 1fil\def\\{\par}{\sc\theauthors}
\par\medskip}

\vglue 0pt plus 0.1fil 

{\small\parskip=0pt\let\newline\\
{\leftskip 0pt plus 1fil\def\\{\par}{\sl\theaddress}\par}
\expandafter\ifx\theemail\relax    
\relax\else\vglue 5pt plus 0.02fil minus 2pt\def\\{\stdspace{\rm 
and}\stdspace} 
\cl{Email:\stdspace\tt\theemail}\fi
\ifx\theurl\relax                  
\relax\else\vglue 5pt plus 0.02fil minus 2pt\def\\{\stdspace{\rm 
and}\stdspace}
\cl{URL:\stdspace\tt\theurl}\fi\par}

\vglue 7pt plus 0.3fil minus 3pt

{\bf Abstract}
\vglue 5pt plus 0.1fil minus 2pt

\theabstract

\vglue 7pt plus 0.3fil minus 3pt

{\bf AMS Classification numbers}\quad Primary:\quad \theprimaryclass

Secondary:\quad \thesecondaryclass

\vglue 5pt plus 0.3fil minus 2pt

{\bf Keywords:}\quad \thekeywords

\vglue 10pt plus 0.5fil minus 5pt

{\small  Proposed: \theproposer\hfill Received: \receiveddate\nl
Seconded: \theseconders\hfill 
\ifx\reviseddate\relax                         
Accepted: \accepteddate                        
\else
Revised: \reviseddate                          
\fi}
\eject
}       

\let\maketitlepage\maketitlep
\let\maketitle\maketitlepage

\font\phead=cmsl9 scaled 950
\font=cmsl9 scaled 1050
\font\pnum=cmbx10 scaled 913
\font\lnum=cmbx10 
\font\pfoot=cmsl9 scaled 950
\font=cmsl9 scaled 1050
\ifplaintex
\headline{\vbox to 0pt{\vskip -4.5mm\line{\small\phead\ifnum
\count0=\startpage ISSN 1364-0380 (on line)
1465-3060 (printed) \hfill {\pnum\folio}\else\ifodd\count0\def\\{ }%
\ifx\theshorttitle\relax\thetitle\else\theshorttitle\fi\hfill{\pnum\folio}
\else\def\\{ and }{\pnum\folio}\hfill\ifx\theshortauthors\relax\theauthors
\else\theshortauthors\fi\fi\fi}\vss}}
\footline{\vbox to 0pt{\vglue 0mm\line{\small\pfoot\ifnum\count0=\startpage
\copyright\ \gtp\hfill\else
\gt, Volume \thevolumenumber\ (\thevolumeyear)\hfill\fi}\vss
}}
\else
\makeatletter
\def\@oddhead{{\small\lhead\ifnum\count0=\startpage ISSN 1364-0380 (on line)
1465-3060 (printed) \hfill {\lnum\number\count0}\else\ifodd\count0
\def\\{ }\ifx\theshorttitle\relax \thetitle \else\theshorttitle\fi\hfill
{\lnum\number\count0}\else\def\\{ and }{\lnum\number\count0}
\hfill\ifx\theshortauthors\relax 
\theauthors\else\theshortauthors\fi\fi\fi}}\def\@evenhead{\@oddhead}
\def\@oddfoot{\small\lfoot\ifnum\count0=\startpage\copyright\ \gtp\hfill\else
\gt, Volume \thevolumenumber\ (\thevolumeyear)\hfill\fi}
\def\@evenfoot{\@oddfoot}
\makeatother
\fi

\newwrite\gtoutfile
\long\gdef\makeheadfile{  
{\def\\{, }\def\s{ }
\immediate\openout\gtoutfile head.xxx
\immediate\write\gtoutfile{Proxy-for: \ifx\theasciiauthors\relax
\theauthors\else\theasciiauthors\fi\s<\ifx\theasciiemail\relax\theemail\else\theasciiemail\fi>}
\immediate\write\gtoutfile{\noexpand\\}
\immediate\write\gtoutfile{Authors: \ifx\theasciiauthors\relax
\theauthors\else\theasciiauthors\fi}
{\def\\{ }\immediate\write\gtoutfile{Title: \ifx\theasciititle\relax
\thetitle\else\theasciititle\fi}}
\immediate\write\gtoutfile{Subj-class: GT or SG or MG etc}
\immediate\write\gtoutfile{MSC-class: \theprimaryclass\ifx\thesecondaryclass\relax\else, \thesecondaryclass\fi}
\immediate\write\gtoutfile{Journal-ref: Geom. Topol. \thevolumenumber
(\thevolumeyear) \startpage-\finishpage}
\immediate\write\gtoutfile{Comments: Published by Geometry and Topology at}
\immediate\write\gtoutfile{\s\s http://www.maths.warwick.ac.uk/gt/GTVol\thevolumenumber/paper\thepapernumber.abs.html}
\immediate\write\gtoutfile{\noexpand\\}
\immediate\write\gtoutfile{}
\ifx\theasciiabstract\relax
\immediate\write\gtoutfile{\theabstract}\else
\immediate\write\gtoutfile{\theasciiabstract}\fi
\immediate\write\gtoutfile{}
\immediate\write\gtoutfile{\noexpand\\}
\immediate\write\gtoutfile{}
\immediate\closeout\gtoutfile}}  

\def\maketitlepage{\maketitlep\makeheadfile}
\let\maketitle\maketitlepage

\lognumber{655}
\received{12 May 2005}
\volumenumber{9}\papernumber{48}\volumeyear{2005}
\pagenumbers{2129}{2158}   
\published{4 November 2005}
\accepted{10 October 2005}
\proposed{Robion Kirby}
\seconded{Cameron Gordon, Wolfgang Lueck}

\usepackage{pinlabel,amsmath,amssymb,psfrag}

\let\Bbb\mathbb

\def\psfraga <#1,#2> #3#4{%
\psfrag {#3}{\smash{\rlap{\kern #1 \raise #2\hbox{#4}}}}}

\def\fref#1{\hyperlink{#1anchor}{\ref*{#1}}}
\def\figref#1{\hyperlink{#1anchor}{Figure~\ref*{#1}}}
\def\anchor#1{\noindent\hypertarget{#1anchor}{\smash{$\phantom{99}$}}}

\theoremstyle{plain}
\newtheorem*{theorem*}{Theorem}
\newtheorem*{lemma*} {Lemma}
\newtheorem*{corollary*} {Corollary}
\newtheorem*{proposition*} {Proposition}
\newtheorem*{conjecture*}{Conjecture}
 \newtheorem{theorem}{Theorem}[section]
\newtheorem{lemma}[theorem]{Lemma}
\newtheorem{corollary}[theorem]{Corollary}
\newtheorem{proposition}[theorem]{Proposition}
\newtheorem{conjecture}[theorem]{Conjecture}
\newtheorem{question}[theorem]{Question}
 {\catcode`@=11\global\let\c@figure=\c@theorem}
\renewcommand{\thefigure}{\thetheorem}
{\catcode`@=11\global\let\c@equation=\c@theorem}

 \theoremstyle{definition}
\newtheorem*{remark}{Remark}
\newtheorem*{definition}{Definition}

\newtheorem*{claim}{Claim}

\theoremstyle{definition}

\def \L {\mathbf{L}}

 \def\eps{\epsilon}

\def\s{\sigma}

\def\Q{\Bbb{Q}}

\def\Z{\Bbb{Z}}

\def\R{\Bbb{R}}

\def\N{\Bbb{N}}
\def\l{\lambda}
\def\Bl{B\ell}
\def\part{\partial}

\def\a{\alpha}

\def\bp{\begin{pmatrix}}
\def\sm{\smallsetminus}
\def\ep{\end{pmatrix}}
\def\bn{\begin{enumerate}}

\def\en{\end{enumerate}}
\def\ba{\begin{array}}
\def\ea{\end{array}}

\def\L{\Lambda}

\def\s{\sigma}
\def\a{\alpha}
\def\b{\beta}
\def\ti{\tilde}
\def\wti{\widetilde}

\def\fr12{\frac{1}{2}}

\def\zt{\Z[t^{\pm 1}]}
\def\into{\hookrightarrow}
\newcommand{\imra}{\looparrowright}
\newcommand{\onto}{\twoheadrightarrow}

\def\ra{\longrightarrow}

\def\sr{\Z\ltimes \Z[1/2]}  

\newcommand{\Arf}{\operatorname{Arf}}
\newcommand{\Hom}{\operatorname{Hom}}

\newcommand{\Ext}{\operatorname{Ext}}
\newcommand{\Wh}{\operatorname{Wh}}

\newcommand{\Ker}{\operatorname{Ker}}
\newcommand{\Rank}{\operatorname{rank}}

\newcommand{\Spin}{\operatorname{Spin}}

\begin{document}
 \title{New topologically slice knots}
\author{Stefan Friedl\\Peter Teichner}
 \address{Department of Mathematics, Rice University\\Houston, TX 77005, USA}
\secondaddress{Department of Mathematics, University of 
California\\Berkeley, CA 94720, USA}
\asciiaddress{Department of Mathematics, Rice University\\Houston, 
TX 77005, US\\and\\Department of Mathematics, University of 
California\\Berkeley, CA 94720, USA}
\gtemail{\mailto{friedl@rice.edu}{\rm\qua
and\qua}\mailto{teichner@math.berkeley.edu}}
\asciiemail{friedl@rice.edu, teichner@math.berkeley.edu}

\begin{abstract}
In the early 1980's Mike Freedman showed that all knots with trivial
Alexander polynomial are topologically slice (with fundamental group
$\Z$). This paper contains the first new examples of topologically
slice knots. In fact, we give a sufficient {\em homological} condition
under which a knot is slice with fundamental group $\sr$. These two
fundamental groups are known to be the only {\em solvable ribbon}
groups. Our homological condition implies that the Alexander
polynomial equals $(t-2)(t^{-1}-2)$ but also contains information
about the metabelian cover of the knot complement (since there are
many non-slice knots with this Alexander polynomial).

\hyperlink{Err}{\it Erratum attached}
\end{abstract}

\asciiabstract{%
In the early 1980's Mike Freedman showed that all knots with trivial
Alexander polynomial are topologically slice (with fundamental group
Z). This paper contains the first new examples of topologically slice
knots.  In fact, we give a sufficient homological condition under
which a knot is slice with fundamental group Z semi-direct product
Z[1/2].  These two fundamental groups are known to be the only 
solvable ribbon groups.  Our homological condition implies that the
Alexander polynomial equals (t-2)(t^{-1}-2) but also contains
information about the metabelian cover of the knot complement (since
there are many non-slice knots with this Alexander polynomial).}


\maketitle

\section{Introduction}
\begin{figure}[ht!]\small\anchor{fig:61-C} 
\begin{center}
\psfraga <0pt,2pt> {C}{$C$}
\includegraphics[scale=0.2]{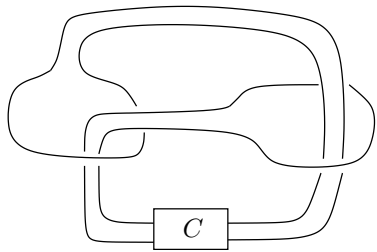}
\end{center}
\caption{A family of ribbon knots} \label{fig:61-C}
\end{figure}

\begin{figure}[ht!]\small\anchor{ribbonsing}
\begin{center}
\includegraphics[scale=0.25]{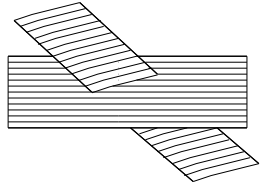}
\caption{Local singularity of a ribbon}
\label{ribbonsing}
\end{center}
\end{figure}

A {\em knot} is an embedding $S^1 \into S^{3}$. We work in the topological category and assume that every embedding is  {\em locally flat}. Note that a smooth embedding is locally flat and that we prefer to {\em add} the adjective ``smooth'' to emphasize this {\em stronger} condition on the embedding, except in our title.
Two knots are called {\em concordant} if there exists an embedding
$S^{1}\times [0,1] \into S^{3}\times [0,1]$ which restricts to the given knots at both ends.
The concordance classes form an abelian group under connected sum, the knot concordance group $\mathcal C$.
A knot is called {\em slice} if it is concordant
to the unknot or, equivalently, if it bounds an embedding of disks $D^{2}\into D^{4}$.
 Predating the 4--dimensional revolution in the early 80's, Casson and Gordon showed that the
epimorphism from $\mathcal C$ onto its high dimension analogue has a nontrivial kernel \cite{CG86}.
There has been much recent progress in understanding how complicated $\mathcal C$ really is. In
\cite{COT03}, \cite{CT04} an infinite sequence of new invariants  was found using non-commutative
Blanchfield forms and their von Neumann signatures.
 On the other hand, many knots are known to be slice, for example the knots in
\figref{fig:61-C}, where the band can be tied into an arbitrary knot $C$.
In fact, there is a large class of
slice knots given as the boundary of {\em ribbons} in $S^3$. These knots are called {\em ribbon
knots}, where a {\em ribbon} is a smooth immersion $D^2\imra S^3$ such that all singularities are
of the type as in \figref{ribbonsing}: They consist of arcs of self--intersection that lie
completely in the interior of one of the two sheets involved. Such singularities can be resolved in
$D^4$ by pushing an open disk around each singular arc slightly away from $\partial D^4$. Thus a ribbon leads to a smooth slice disk in $D^4$, the so called {\em ribbon disk}.
 It is a fascinating open problem whether every smoothly slice knot is ribbon. One distinctive feature
of a ribbon knot is that the inclusion map induces an {\em epimorphism} of the knot group onto the
{\em ribbon group}, the fundamental group of the complement in $D^4$ of the ribbon disk. There is a
simple criterion for a given group to be ribbon in terms of certain presentations, see
Theorem~\ref{thm:ribboncrit}. For general slice complements the inclusion map does {\em not} induce an epimorphism,  and slice disks with the ontoness property are called {\em homotopically ribbon}, or {\em h--ribbon} for short.

In the topological category, Freedman proved that any knot with trivial  Alexander polynomial is
slice \cite{FQ90}, see also \cite{GT04} for a more direct construction. Using gauge theory,
Gompf showed that some of these knots are {\em not} smoothly slice \cite{G86}. The easiest such
knot is probably the Whitehead double of the trefoil knot, already exhibiting the subtle
difference between smooth and topological 4--manifolds. In an amazing turn of events, Rasmussen
very recently gave the first purely combinatorial proof for the fact that this knot is not
smoothly slice. He constructed a concordance invariant from Khovanov homology \cite{R04} with
beautiful properties. In particular, the arguments of \cite{L04} for showing that (iterated)
Whitehead doubles of the trefoil are not smoothly slice can be adapted to this setting from the
concordance invariant of Ozsv\'ath and Szab\'o \cite{OS03}. Livingston uses the following
argument, going back to at least Rudolph: One can exhibit some Whitehead doubles as separating
curves on minimal Seifert surfaces of certain torus knots. Since the Rasmussen and
Ozsv\'ath--Szab\'o invariants detect the minimal genus of torus knots, it follows that such
separating curves cannot be smoothly slice.

In this note we provide the first new class of slice knots since Freedman's construction,
using
his theorem \cite{F82} that solvable groups are good (for topological surgery).
 Our main result is the following theorem. Let
\[
SR:= \langle a, c \mid aca^{-1}=c^2 \rangle \cong \sr.
\]
Here the generator $a$ of $\Z$ acts on the normal subgroup $\Z[1/2]$ via multiplication by $2$. It is known, cf Lemma~\ref{lem:LSR-groups}, that $SR$ and $\Z$ are the only {\em solvable ribbon} groups, hence the name.  In geometric group theory, this group is also known as the {\em Baumslag--Solitar group} $B(1,2)$.

\begin{figure}[ht!]\small\anchor{example2}
\begin{center}
\psfrag {A}{$A$}
\begin{tabular}{ccc}
\includegraphics[scale=0.22]{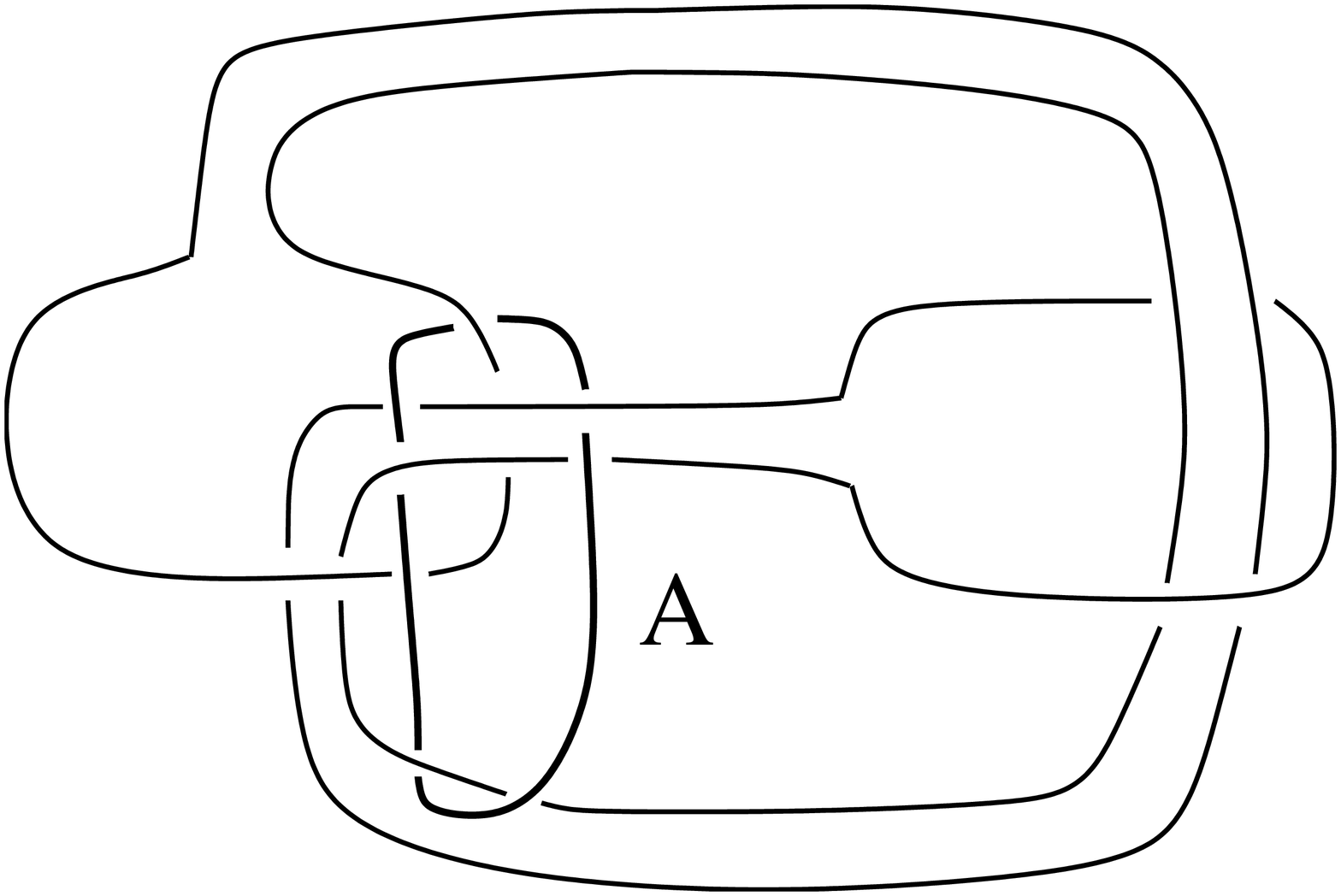}&\hspace{1cm}&
\includegraphics[scale=0.275]{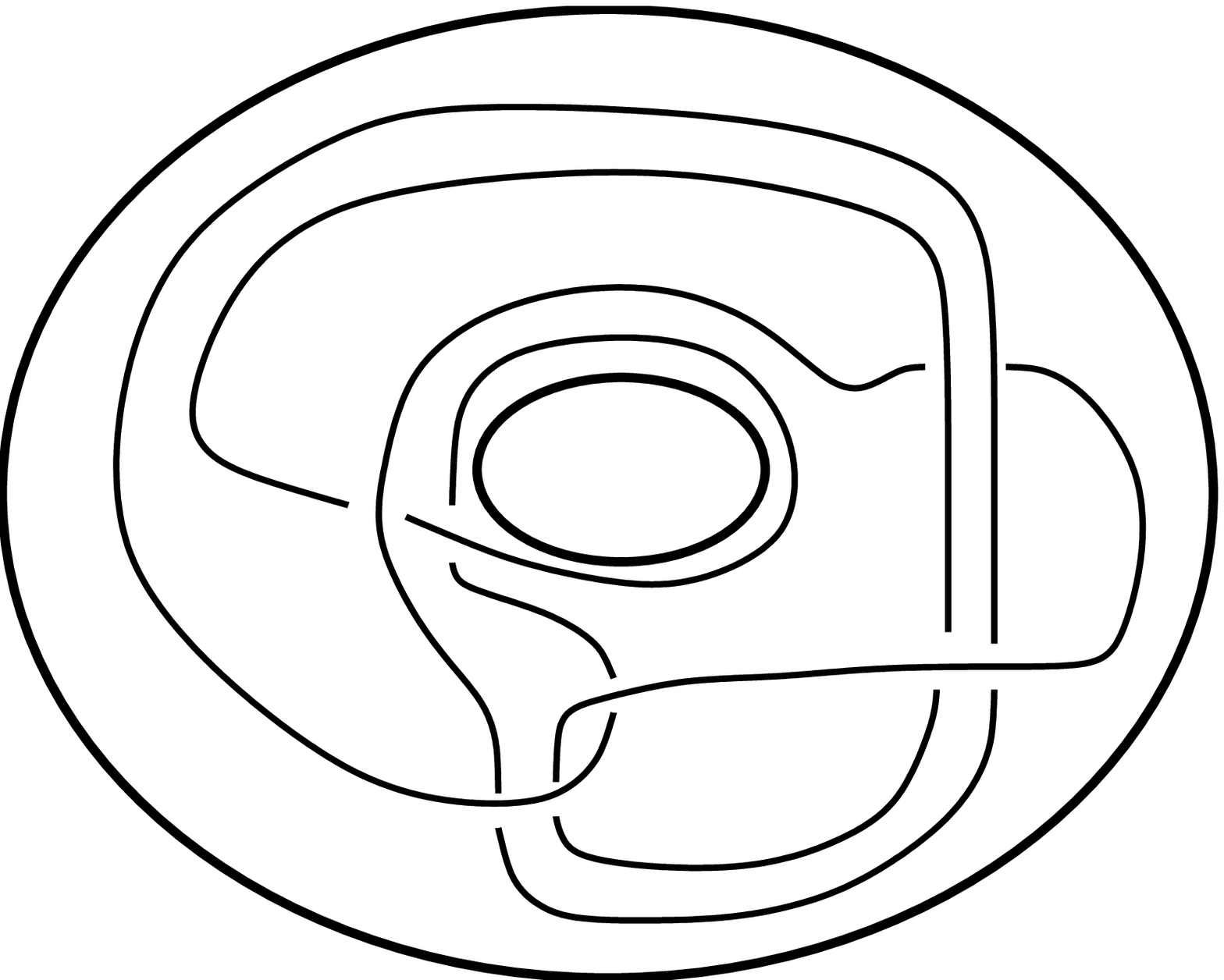}
\end{tabular}
\caption{The knot $6_1$ in the solid torus $S^3\sm A$}
\label{example2}
\end{center}
\end{figure}

\begin{theorem}
\label{thm:main} Let $K$ be a knot and denote by $M_K$ the $0$--surgery on $K$. If
there is an epimorphism $\pi_1(M_K) \onto G$, $G=SR$ or $\Z$, such that \begin{equation*}\tag{Ext}
\Ext_{\Z[G]}^1(H_1(M_K;\Z[G]),\Z[G])=0 \end{equation*} then $K$ is (topologically) slice. In fact,
$K$ is h--ribbon with group $G$ if and only if this Ext--condition holds for some epimorphism
$\pi_1(M_K) \onto G$.
\end{theorem}
 The case $G=\Z$ is actually just a reformulation of Freedman's theorem because the condition (Ext)
is then equivalent to $\Delta_K(t)=1$. For $G=SR$, we shall show in Corollary \ref{cor:arf}
 that this condition implies \[
\Delta_K(t) = (t-2)(t^{-1}-2). \]
There are well known knots with this Alexander polynomial that
are not slice (cf Section~\ref{sec:61}) so the h--ribbon question with group $SR$ is more subtle
than for $\Z$.
Our result complements work of Tim Cochran and Taehee Kim \cite{CK05}. They show that if the degree of the Alexander polynomial
is greater than two, than the homology of solvable covers can not determine whether a given knot is slice or not.
 \begin{remark} We will show in Lemma~\ref{lem:Bl=0} that the somewhat awkward (Ext)
condition can be replaced by the condition that a non-commutative Blanchfield pairing \[ \Bl\co
H_1(M_K;\Z[G])\times H_1(M_K;\Z[G])\ra Q(G)/\Z[G]
\]
vanishes, where $Q(G)$ denotes the Ore localization of $\Z[G]$.
\end{remark}
 It is surprisingly easy to construct many knots that satisfy all conditions of our
Theorem~\ref{thm:main}. For example, all knots  in \figref{fig:61-C} do.
The easiest one is the knot $6_1$, which is isotopic to the case where the band $C$ is the unknot. However, all knots in
\figref{fig:61-C} are obviously ribbon and hence we need to work
harder to get new h--ribbon knots. One trick is the following satellite construction
that can be considered as an analogue of Whitehead doubling (a nice way to
obtain knots with trivial Alexander polynomial).
We give a much more general result in Section~\ref{sec:satellite}.
 \begin{theorem} \label{thm:Wh}
Start with the knot $6_1$, drawn in a solid torus as on the right hand side of \figref{example2}. Tie an arbitrary knot into that solid torus to obtain a satellite $K$ of $6_1$. Then $K$ satisfies the assumptions of Theorem~\ref{thm:main} and is therefore h--ribbon.
\end{theorem}
 On the left hand side of \figref{example2} we picked an axis $A$ which is unknotted in $S^3$ (so that the complement of $A$ will be a solid torus) and which punctures all ribbon disks that we could see. Then we redrew the picture in a way that $A$ becomes the meridian to the visible solid torus.
 We conjecture that Theorem~\ref{thm:Wh} gives new examples of knots
which are h--ribbon but not smoothly slice. The simplest candidate is the satellite knot of the trefoil knot. It has a knot diagram with 93 crossings and we are unable to compute
Rasmussen's invariant \cite{R04}, the complexity of which is exponential in the
number of crossings.
For the general ribbon case we propose the following generalization of Theorem~\ref{thm:main}.
 \begin{conjecture} \label{conj:ribbon}
Let $G$ be a ribbon group for which topological surgery works.
A knot $K$ is h--ribbon with group $G$ if and only if there exists an epimorphism $\varphi\co \pi_1(M_K)\onto G$ such that
the condition (Ext) from Theorem~\ref{thm:main} holds.
\end{conjecture}
Note that the converse of this conjecture is not completely straight forward either. We use a Blanchfield form to prove the converse for the groups $G=SR$ and $\Z$.
 \begin{question}
Are the fundamental groups of complements of h--ribbons ribbon groups?
\end{question}
The phrase ``surgery works for $G$'' means here that the (reduced) surgery sequence
\begin{equation}\tag{$\mathcal S$}
\mathcal S_{TOP}^h(X,M)\ra  \widetilde {\mathcal N}_{TOP}(X,M)\ra \widetilde L_4^h(\Z[G])
\end{equation}
is exact for all Poincar\'e pairs $(X,M)$ with $\pi_1(X)=G$, see \cite[chapter~11.3]{FQ90} for more information. If $G=SR$ or $\Z$ then this sequence is exact by Freedman's disk embedding theorem (cf \cite{F82} and \cite{FT95}) for solvable groups. The general case is still open and since most other ribbon groups contain free groups we need a very strong form of this result.
It is logically possible that topological surgery works for a given fundamental group but the disk embedding theorem (and hence the s--cobordism theorem) fail for this group.

The ``if''--direction of Conjecture~\ref{conj:ribbon} would follow from the following purely homological conjecture.
We say  that $(M_K,\varphi)$ satisfies the  {\em Poincar\'e duality} condition
if the induced inclusion $M_K\into K(G,1)$ is a finite 4--dimensional Poincar\'e pair.
\begin{conjecture} \label{conj:PD}
Let $K$ be a knot with an epimorphism $\pi_1(M_K)\onto G$ onto a ribbon group such that the
condition (Ext) from Theorem~\ref{thm:main} holds. Then $(M_K,\varphi)$ satisfies the Poincar\'e
duality condition. Moreover, $ \widetilde L_4^h(\Z[G])=0$.
\end{conjecture}
We give supporting evidence for this  algebraic conjecture at
the end of Section~\ref{sec:ribboncomp}. For $G=\Z$ it is true out of easy reasons and in Lemmas~\ref{lem:pcext},  \ref{lem:LSR} and Lemma~\ref{lem:LSR4} we prove that the conjecture is true for $G=SR$. The main reason why this works is the fact that the group ring $\Z[SR]$ is
an Ore-domain and hence has an ordinary (skew) quotient field, see Section~\ref{sec:SR}. We use
various lemmas and ideas from \cite{COT03}. The first statement of Theorem~\ref{thm:main} is then a consequence
of the following general result.
 \begin{theorem} \label{thm:mainG}
Let $G$ be a finitely presented group for which topological sur\-gery works and with $H_1(G) \cong \Z, H_2(G)=0$ and $ \widetilde L_4^h(\Z[G])=0$.
A knot $K$ is h--ribbon with group $G$ if there is an epimorphism $\varphi\co \pi_1(M_K)\onto G$  such that $(M_K,\varphi)$ satisfies the Poincar\'e duality condition.
In particular, Conjecture~\ref{conj:ribbon} follows from Conjecture~\ref{conj:PD}.
\end{theorem}
 Theorem~\ref{thm:mainG} comes from the fact that given a Poincar\'{e} pair, one can attempt to use classical surgery theory to find a $4$--manifold $W$ of type $K(G,1)$. We show that in the situation above, this approach successfully produces a h--ribbon complement that is a
$K(G,1)$.

 The paper is organized as follows. In Section \ref{sec:ribboncomp} we recall several known facts about ribbon groups.
In Section \ref{sec:proof}  we explain how surgery can be used to find h--ribbons. We accumulate several restrictions
on the group $G$, and in Section \ref{sec:SR} we show that the group $SR$ satisfies all of them.
In Section  \ref{sec:BL} we will show that for our groups the (Ext)--condition can be replaced
by a vanishing condition on a non--commutative Blanchfield pairing.
In Section  \ref{sec:satellite} we recall the satellite construction
which we use to give examples of topologically slice knots in
Section \ref{sec:61}.

\medskip{\bf Acknowledgments}\qua We thank Chuck Livingston and the referee for their
valuable comments on the first version of this paper. We also wish to thank Jerry Levine and
Andrew Ranicki for helpful discussions regarding certain technical aspects in this paper.

The first author was partially supported by RTN Network HPRN-CT-2002-00287: Algebraic K--Theory,
Linear Algebraic Groups and Related Structures. The second author is partially supported by NSF
Grant DMS-0453957.

\section{Ribbon disk complements} \label{sec:ribboncomp}
\begin{figure}[ht!]\small\anchor{figkirby} 
\begin{center}
\psfrag {a}{$a$}
\psfrag {b}{$b$}
\psfraga <0pt,2pt> {C}{$C$}
\includegraphics[scale=0.22]{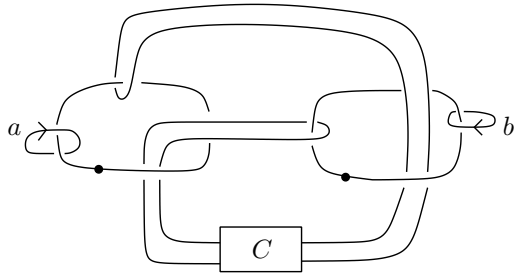}
\end{center}
\caption{Handle decomposition of $N_D$, with generators $a,b$ for $\pi_1(N_D)$.}
\label{figkirby}
\end{figure}

 \begin{definition}
A group $G$ is called {\em ribbon} if there exists a ribbon disk
$D \into D^4$, as explained in the introduction, with $\pi_1(D^4\sm D)\cong G$.
\end{definition}
 The deficiency of a presentation of a group is the number  of generators minus the number of
relations. The deficiency $\mbox{def}(G)$ of a finitely presented group $G$ is the maximum of the
deficiencies of all presentations. The following theorem is well-known and uses the fact that any
ribbon in 3--space can be obtained as follows: Start with an $s$--component unlink and add $(s-1)$
bands to produce a knot, where the bands can be arbitrarily twisted and linked. In particular, the
bands may hit the disks bounding the unlink which introduces our ribbon singularities. This
description actually gives an embedded disk $D\into D^4$ with a Morse function $D^4\to [0,1]$ which
has $s$ minima and $(s-1)$ saddles when restricted to $D$. This implies that the complement
$N_D:=D^4 \sm \nu D$, $\nu D$ an open tubular neighborhood of $D$, has a handle decomposition with one $0$--handle, $s$ 1--handles and
$(s-1)$ 2--handles, giving the desired group presentation.
 \begin{theorem} \label{thm:ribboncrit}
A group $G$ is ribbon if and only if
$H_1(G)=\Z$ and $G$ has a Wirtinger presentation of deficiency one, ie, $G$ has a presentation
of the form
\[ G=\langle g_1,\dots,g_s \mid r_1,\dots,r_{s-1} \rangle, \]
where $r_i=g_{h_i}g_{l_i}^{\eps_i}g_{k_i}^{-1}g_{l_i}^{-\eps_i}$ for some $h_i,k_i,l_i \in \{1,\dots,s\}$
and $\eps_i\in \{-1,1\}$.
 \end{theorem} Using the above theorem and the  Wirtinger presentation, one sees that knot groups
are ribbon. As pointed out by John Stallings, a direct geometric way to see the corresponding
ribbon disk is as follows. Represent the knot by an arc ending in a plane $P\subset\R^3$. A
rotation by 180 degrees in $\R^4$, with fixed plane $P$, sweeps out a disk in $D^4$ which has the
same local minima and saddles as the original knotted arc (when projected orthogonally to the fixed
plane) and there are no maxima. Applying this construction to a knot $K$, one obtains a ribbon disk
with fundamental group $\pi_1(S^3 \sm K)$ and boundary $K \# (-K)$. We recover by this argument the
reason why the mirror image $-K$ is the inverse of $K$ in the concordance group.
 \begin{lemma} \label{lem:LSR-groups}
The only solvable ribbon groups are $\Z$ and $SR=\sr$.
\end{lemma}
 \begin{proof}
Wilson \cite{W96} shows that every solvable group of deficiency one
is isomorphic to $G_k:=\Z \ltimes \Z[t,t^{-1}]/(tk-1)$ for some $k$.
A direct computation shows that $G_k/[G_k,G_k]=\Z$ if and only if
$k=0,2$. These are exactly the two groups in the statement. Below we show that $SR$ is indeed a ribbon group.
\end{proof}
 We now show that $SR$ does in fact appear as the fundamental group of a ribbon disk complement. Let
$D$ be the ribbon disk given by \figref{fig:61-C}. The corresponding ribbon disk complement
$N_D$ is given by the handle diagram of \figref{figkirby}, with two (dotted) 1--handles and one
2--handle. This is the case $s=2$ in the discussion above Theorem~\ref{thm:ribboncrit}, see also
\cite[p.~213]{GS99}. 
 Note that with $c:=ab^{-1}$ one gets \[\pi_1(D^4 \sm D)=\langle a,b\mid
a^{-1}ba^{-1}bab^{-1}\rangle = \langle a,c \mid aca^{-1}=c^2\rangle =SR. \] which is in fact
independent of the knot $C$ that was tied into the band in \figref{fig:61-C}.
 \begin{remark}
Writing $R(C)$ for the above ribbon knot, then $R(C)$ is a {\em satellite knot}
with companion $C$ and orbit $R:=R(\mbox{trivial knot})$.
The knot $R$ is the knot $6_1$ of the standard knot table
and we shall construct many more examples in this manner in Section~\ref{sec:satellite}.
\end{remark}
 \begin{remark} Since any ribbon group $\pi_1(N_D)$ has a presentation of deficiency
$\Rank(H_1(N_D))=1$ it is in particular an E--group (cf \cite[p.~324]{S74}), hence
$\pi_1(N_D)^{(1)}$ is an E--group (cf \cite[p.~302]{S74}). It follows from a result of Roushon
\cite[Corollary~4.6]{R03} in the case that $\Rank_{\Z}(\pi_1(N_D)^{(1)}/\pi_1(N_D)^{(2)})\geq 2$ that the
derived series of $\pi_1(N_D)$ never stabilizes, ie, $\pi_1(N_D)^{(i)} \ne \pi_1(N_D)^{(i+1)}$
for all $i\in \N$ (cf also \cite{C04}). This means that with few exceptions ribbon groups are
neither solvable nor does their derived series stabilize. \end{remark}

\subsection*{The ribbon group conjecture} \label{sec:conj}
Alexander duality implies that ribbon groups satisfy $H_1(G) \cong \Z$ and $H_2(G)=0$. We say that a group $G$ is {\em aspherical} if $K(G,1)$ is a 2--complex.
For a ribbon group  $G=\pi_1(D^4\sm D)$ this is the case if $\pi_2(D^4\sm D)=0$.
It is conjectured that ribbon groups are aspherical. Note that this is in turn
a special case of the Whitehead conjecture.
This conjecture is known to be true for all knot groups
since knot complements are aspherical by the sphere theorem. It is also known for
all locally indicable ribbon
groups (cf \cite{H82}), which in particular includes all ribbon groups with
$G^{(\a)}=\{e\}$ for some ordinal $\a$ (cf also \cite{S74}). For some more examples
cf \cite{H85}.
 In particular, we conjecture that a ribbon group $G$ has the property that \[ H_3(G)=0 \mbox{ and }
\Ext_{\Z[G]}^i(\Z,\Z[G])=0 \mbox{ for }i>2. \] This conjecture clearly holds for aspherical ribbon
groups. Furthermore we conjecture that for a ribbon group $G$ and an epimorphism
$\varphi\co \pi_1M_K\onto G$, we always have \[ \Hom_{\Z[G]}(H_i(M_K;\Z[G]),\Z[G])=0 \mbox{ for }
i=1,2. \] The relevance of these properties comes from Lemma~\ref{lem:pcext} and
Theorem~\ref{thm:mainG}. We will prove all these properties for $G=SR$ in Lemma~\ref{lem:LSR}.

\vspace{-0.2cm}
\section{Proof of the Main Theorem~\ref{thm:main}} \label{sec:proof}
\vspace{-0.2cm}
 Let $K$ be a knot in $S^{3}$ and denote by $M_K$ the result
of zero framed surgery along $K$. The following is a well known slice criterion,
see eg, \cite{COT03}. For the ``only if'' direction one takes $W$ to be
the complement in $D^4$ of a (thickened up) slice disk. For the ``if'' direction
one uses Freedman's solution of the topological Poincar\'{e} conjecture in dimension~4,
in order to recognize as $D^4$ the $4$--manifold $W$ union a 2--handle along a meridian for $K$ in $M_K$.
 \begin{proposition} \label{prop:topslice}
A knot $K$ is slice if and only if $M_K^{3}$ bounds a $4$--manifold $W^{4}$ with
\bn
\item $\pi_1(W)$ is normally generated by the image of a meridian for $K$,
\item $H_1(W)\cong \Z$,
\item $H_2(W)=0$,
\en
Furthermore, $K$ is h--ribbon with group $\pi_1(W)$ if and only if $(1)$ is replaced by
\bn
\item[\rm(1-h)\hspace{-3pt}] $\pi_1(M_K) \to \pi_1(W)$ is surjective.
\en
\end{proposition}
 Note that the first two conditions can always be satisfied, even with $\pi_1(W)\cong\Z$. However, to satisfy condition (3) as well, it is often necessary to make the fundamental group of $W$ more complicated. In fact, if (1)--(3) are satisfied for $\pi_1(W)\cong\Z$ then $K$ has vanishing Alexander polynomial. Thus the main problem is to find a candidate for the fundamental group of $W$.

 \subsection{Outline of the construction of W} \label{sec:outline}
To prove Theorem~\ref{thm:main}, we take our candidate for the fundamental group to be $G=SR$ or $\Z$. By assumption we have given an epimorphism $\varphi\co  \pi_1M_K \onto G$ and we would like to find a $K(G,1)$--manifold $W$ with boundary $M_K$ such that the inclusion induces $\varphi$. Since any ribbon group
satisfies $H_1(G)=\Z$ and $H_2(G)=0$, such a manifold $W$ would fulfill the conditions (1-h), (2) and (3) above. This means that we are in a classical surgery situation and we can follow the steps taken by many before us. Below, the steps are labelled by the section numbers where they will be worked out.

\begin{itemize}
\item[(\ref{sec:PC})]
 First we seek conditions such that $(K(G,1), M_K)$ is a Poincar\'{e} pair. It turns out that exactly condition (Ext) from Theorem~\ref{thm:main} arises.
\item[(\ref{sec:degree 1})]
 Secondly, we check that the Spivak normal bundle of this Poincar\'{e} pair has a linear reduction, ie, that there is a degree~1 normal map from a (smooth) manifold pair $(N,M_K)\to (K(G,1), M_K)$. In fact, we shall see that there is a {\em unique} normal cobordism class of such maps if we require, as we shall, that the signature of $N$ is zero.
\item[(\ref{sec:surgery})]
 This normal map has surgery obstruction in the reduced L--group\break $\widetilde L^h_4(\Z[G])$. Since $\pi_2K(G,1)=0$ this quadratic form is nothing but  the intersection form  on $\pi_2N$. We will verify in Lemma~\ref{lem:LSR4} that our L--group vanishes and hence the intersection form  on $\pi_2N$ is (stably) hyperbolic.
\item[(\ref{sec:final})]
Finally, since by assumption the surgery sequence $(\mathcal S)$ for our fundamental groups is exact (in the topological category), one gets a topological manifold $W$ together with a homotopy equivalence $(W,M_K)\to (K(G,1), M_K)$. Hence we found our h--ribbon by Proposition~\ref{prop:topslice}.
\end{itemize}
In the following sections, we go through these steps one by one. We will {\em not} assume from the beginning that $G=SR$ or $\Z$ but we shall add conditions on $G$ (and $\varphi$) as we go, in turn proving Theorem~\ref{thm:mainG}. It will turn out that some of these conditions are known for all ribbon groups and some are conjectural. For the groups $SR$ and $\Z$ we will prove all conditions in Section~\ref{sec:SR}.

\subsection{The Poincar\'e duality condition} \label{sec:PC}
We will use the following conventions. If $\ti{X} \to X$ is the universal cover then $C_*(\ti{X})$
has a right $\Z[\pi_1(X)]$--module structure. Given a left $\Z[\pi_1(X)]$--module $P$ we can
consider $H_*(X;P):=H_*(C_*(\ti{X})\otimes_{\Z[\pi_1(X)]}P)$. Using the usual involution on group
rings given by $\bar{g}=g^{-1}$ we can turn left $\Z[G]$--modules into right $\Z[G]$--modules. We
define $H^*(X;P):=H_*(\mbox{Hom}_{\Z[\pi_1(X)]}(C_*(\ti{X}),P))$.

 \begin{definition}\cite{W99}\qua
Let $Y\subset X$ be CW--complexes, write $G:=\pi_1(X)$. Then $(X,Y)$
is called a   Poincar\'e pair of dimension $n$ if there
exists $[X]\in H_n(X)$, such that
\[
\ba{rrcl}\tag{$\cap$}
\cap [X]\co & H^i(X,Y;\Z[G])&\ra &H_{n-i}(X;\Z[G]) \\
\cap \partial [X]\co & H^i(Y;\Z[G])&\ra &H_{n-i-1}(Y;\Z[G])
\ea
\]
are isomorphisms of $\Z[G]$--right modules. $[X]$ is called the fundamental homology class
of $(X,Y)$.
\end{definition}
 Note that if $(X,Y)$ is a Poincar\'e complex, then from the long exact homology and cohomology sequences
 it follows that
\[ \cap [X]\co  H^i(X;\Z[G])\to H_{n-i}(X,Y;\Z[G])\] is an isomorphism as well. In our case, $Y$ will
be the 3--manifold $M_K$ so that the second condition on $\partial [X]$ will be satisfied.
 Recall that we have an inclusion $\varphi\co M_K \into K(G,1)$ that induces an epimorphism of
fundamental groups. For brevity we set
\[ M:= M_K \quad \text{ and } \quad X:=K(G,1). \]
Recall
that if $(X,M)$ is a  finite 4--dimensional Poincar\'e pair, we say that $(M,\varphi)$ satisfies the
{\em Poincar\'e duality} condition.
 \begin{lemma}\label{lem:pc}
Assume that $G$ satisfies $H_3(G)=0$ and $H^i(G;\Z[G])=0$ for $i>2$.
Then $(M,\varphi)$ satisfies the Poincar\'e duality condition if and only if
\[ \varphi^*\co  H^i(G;\Z[G])\ra H^i(M;\Z[G])
\]
is an isomorphism for $i=1,2$.
\end{lemma}
 \begin{remark} In Section~\ref{sec:SR} we will show that the above assumptions are satisfied for
$G=SR$. In this case $H^2(G;\Z[G])\ne 0$ (cf Lemma  \ref{lem:hig}) which shows that the Poincar\'e duality condition in
general can not be simplified to $H^2(M_K;\Z[G])=0$. Note also that the lemma implies that
$H_i(G)=0$ for all $i\geq 3$. \end{remark}

\begin{proof}
We first show the ``only if'' direction. In that case $H^{4-i}(X,M;\Z[G])\cong H_{i}(X;\Z[G])$. But $H_i(X;\Z[G])=0$ for $i\neq 0$ since  this is the homology of the $G$-cover of $X$ which is contractible since $X=K(G,1)$. The claim now follows from the long exact cohomology sequence with $\Z[G]$--coefficients.
 We now turn to the proof of the ``if'' direction.
By the long exact homology sequence of $(X,M)$, with $\Z$--coefficients, and the vanishing of
$H_3(G)=H_3(X)$ we can choose a class $[X]\in H_4(X,M)$ that maps to the fundamental class $[M]
\in H_3(M)$ under the boundary map. Cap product with this class induces maps as in $(\cap)$
above and we need to check that they are isomorphisms. These maps are the left (or right) most
vertical arrows in the commutative diagram of long exact sequences 
(with $\Z[G]$--coefficients
understood):
\[ \ba{rccccccc}
\hspace{0.0cm} \to \hspace{0.0cm}& H^i(X,M)&\hspace{0.0cm}\to\hspace{0.0cm} &
H^i(X)&\hspace{0.0cm}\to \hspace{0.0cm}& H^i(M)&\hspace{0.0cm}\to \hspace{0.0cm}& H^{i+1}(X,M)
\\
[0.2cm] &\downarrow \hbox{\scriptsize $\cap (-1)^i[X]$} & &\downarrow \hbox{\scriptsize $\cap (-1)^i[X]$}
& &\downarrow \hbox{\scriptsize $\cap (-1)^{i+1}[M]$} &&\downarrow 
\hbox{\scriptsize $\cap (-1)^{i+1}[X]$}
\\
[0.2cm]  \hspace{0.0cm}\to\hspace{0.0cm} & H_{4-i}(X)&\hspace{0.0cm}\to\hspace{0.0cm}&
H_{4-i}(X,M) &\hspace{0.0cm}\to \hspace{0.0cm} & H_{3-i}(M)&\hspace{0.0cm}\to\hspace{0.0cm} &
H_{4-(i+1)}(X)
\ea \]
Recall that $H_i(X;\Z[G])=0$ for $i\neq 0$. Therefore, we need to show in particular that
$H^i(X,M;\Z[G])=0$ for $i\neq 4$. Note that
\[
H_{4-i}(X,M;\Z[G]) \cong H_{3-i}(M;\Z[G]) \quad \text{ for } i\neq 3.
\]
Since the maps $\cap [M]$ are isomorphisms it suffices to show that for all $i\neq 3$ the maps
\[
\varphi^*\co H^i(X;\Z[G]) \ra H^i(M;\Z[G])
\]
are isomorphisms. For $i=1,2$ this is our assumption so it suffices to discuss the other cases. For $i=0$ both groups are zero since $G$ has infinite order. For $i>3$ again both groups are zero, for $X$ this is our assumption and for $M$ it follows from $\dim M=3$.
 Finally, we need to discuss the special case $i=3$ in the diagram above. Then
\[
H_0(M;\Z[G]) \cong \Z \cong H_0(X;\Z[G])
\]
and by assumption $H^i(X;\Z[G])=0$ for $i=3,4$. It follows that the boundary map $H^3(M;\Z[G]) \to H^4(X,M;\Z[G])$ is an isomorphism and by commutativity of the diagram the last map in question
\[
\cap [X]\co  H^4(X,M;\Z[G]) \ra H_0(X;\Z[G])
\]
is also an isomorphism.
\end{proof}

\begin{lemma}\label{lem:pcext}
Let $G$ satisfy $H_3(G)=0$ and $H^i(G;\Z[G])=0$ for $i>2$. Furthermore, assume that $\Hom_{\Z[G]}(H_i(M;\Z[G]),\Z[G])=0$ for $i=1,2$.
Then $(M,\varphi)$ satisfies the Poincar\'e duality condition if and only if
the (Ext) condition from Theorem~\ref{thm:main} holds:
$\Ext_{\Z[G]}^1(H_1(M;\Z[G]),\Z[G])=0$.
\end{lemma}

 \begin{proof}
Denote the universal cover of $M$ by $\wti{M}$.  Write $\pi:=\pi_1(M)$.
Note that we have a chain isomorphism of right $\Z[G]$--module complexes given by
\[\ba{rcl} \Hom_{\Z[\pi]}(C_*(\wti{M}),\Z[G]) &\to & \Hom_{\Z[G]}(C_*(\wti{M})\otimes_{\Z[\pi]}\Z[G],\Z[G])\\
\phi &\mapsto& (c\otimes f \mapsto \overline{\phi(c)}f ).\ea \] Therefore $H^{*}(M;\Z[G])\cong
H_*(\Hom_{\Z[G]}(C_*(\wti{M})\otimes_{\Z[\pi]}\Z[G],\Z[G]))$. We now apply the  universal
coefficient spectral sequence (UCSS). This has an $E_2$--term
\[
E_2^{p,q}=\Ext_{\Z[G]}^p(H_q(M;\Z[G]),\Z[G]),
\]
differentials $d_r$ of degree $(r,1-r)$ and converges to
\[
H^{p+q}(M;\Z[G])\cong
H_{p+q}(\Hom_{\Z[G]}(C_*(\wti{M})\otimes_{\Z[\pi]}\Z[G],\Z[G])),
\]
 cf \cite[Theorem~2.7]{L77} for
more details. We use that the edge homomorphism at $q=0$ of the spectral sequence is the map
\[
\varphi^* \co  \Ext_{\Z[G]}^p(H_0(M;\Z[G]),\Z[G]) \cong H^p(G;\Z[G]) \ra H^p(M;\Z[G]).
\]
By Lemma~\ref{lem:pc} we only need to verify that this is an isomorphism for $p=1,2$ if and only if condition (Ext) holds. For $p=1$ this follows immediately from our assumption
$\Hom_{\Z[G]}(H_1(M;\Z[G]),\Z[G])=0$. For $p=2$ we first observe that since $H^3(G;\Z[G])=0$ there are no possible $d_2$--differentials into the $E_2^{1,1}$ spot.
Since $\Hom_{\Z[G]}(H_2(M;\Z[G]),\Z[G])=0$, we get a short exact sequence
\[
0\ra \Ext_{\Z[G]}^1(H_1(M;\Z[G]),\Z[G]) \ra H^2(G;\Z[G])\overset{\varphi^*}{\ra} H^2(M;\Z[G])\ra 0
\]
showing that the condition (Ext) holds if and only if $\varphi^*$ is an isomorphism.
\end{proof}
 Let $G$ be any group with $H_1(G)=\Z$. Then we define $\Delta_G(t)\in \zt$ to be the order of
the $\Z[t^{\pm 1}]$--module $H_1(G,\Z[t^{\pm 1}])\cong H_1(G^{(1)},\Z)$.
Note that $\Delta_G$ is well-defined up to multiplication by a unit in $\zt$ and up taking the natural involution
$t\mapsto t^{-1}$.
For example if $G=SR$, then $\Delta_G(t)=t-2$.
 \begin{corollary}\label{cor:arf}
Let $G$ and $M$ satisfy all the conditions from Lemma~\ref{lem:pcext}, including the (Ext) condition. Then
\[ \Delta_K(t)=\Delta_G(t)\Delta_G(t^{-1}).\]
In particular $\Arf(K)=0$.
\end{corollary}
 \begin{proof}
By Lemma \ref{lem:pcext} $(K(G,1),M_K)$ is a Poincar\'e pair. Since
\[
H_i(M_K;\Z)\xrightarrow{\cong} H_i(K(G,1);\Z)
\]
for all $i$ it follows from a well-known argument
that
\[
0=H_i(M_K;\Q(t))\xrightarrow{\cong} H_i(K(G,1);\Q(t))
\]
for all $i$. In particular $H_i(K(G,1);\zt)$ is $\zt$--torsion for all $i$.
The corollary now follows from the long exact sequence of $(K(G,1),M_K)$ with $\zt$--coefficients, duality,
and the fact that the alternating product of orders of a long exact sequence equals one (cf \cite[Lemma~5, p.~76]{L67}).
 The Arf invariant $\Arf(K)\in \Z/2$ of a knot $K$ has the well-known property that
it equals zero if $\Delta_K(-1)\equiv \pm 1\mod 8$, and 1 otherwise.
 \end{proof}
 We shall now work through the remaining steps outlined in Section~\ref{sec:outline}.
We have a Poincar\'{e} pair $(X,M)$ where $X=K(G,1)$ and $M=M_K$.
Clearly $\pi_1(X)$ is normally generated by the image of a meridian for $K$ since $\varphi$ is surjective and since any meridian normally generates
$\pi_1(M_K)$.
Furthermore $H_1(X)\cong H_1(G)\cong \Z$ and $H_2(X)=H_2(G)=0$.
It follows from the homology exact sequence of the pair $(X,M)$, Poincar\'e duality, and the fact that $\varphi$ is an isomorphism on $H_1$ that $H_3(X)=0$.
\subsection{The degree 1 normal map} \label{sec:degree 1}
We show that the stable Spivak normal bundle of $(X,M)$ has a linear reduction. It seems well
known to experts that this is always true in the oriented case in 4 dimensions but we could not
find an explicit reference. So here is a rather ad hoc argument: Since $M$ is orientable, its
tangent bundle is trivial and so is its stable normal bundle. We claim that the stable Spivak
normal bundle of $X$ is also trivial: Let $BS$ be the classifying space of such (oriented) bundles.
Its homotopy groups are the (shifted) stable homotopy groups of spheres, starting with
\[ \pi_2BS
\cong\Z/2\cong\pi_3BS\quad\text{ and }\quad \pi_4BS \cong\Z/24.
\]
It follows that for all $i$ we
have $H^i(X;\pi_iBS)=0$ which shows that the classifying map $X\to BS$ of the stable Spivak normal
bundle is trivial.  Finally, it is clear that, up to homotopy, this trivial map can be lifted (trivially) through the fibration
$BSO\to BS$, having the trivial map on the boundary $M$.
From the usual transversality theory we obtain a degree~1 normal map from a manifold pair
\[
f\co (N,M)\to (X, M)
\]
that is the identity on the boundary. Note that since
the stable normal bundle of $X$ is trivial, so is that of $N$. In particular, $N$ is spin. Another approach to finding $f$ is to prove that $(M,\varphi)$ represents the zero element in the spin-bordism group $\Omega^{\Spin}_4(X)$. The fact that this element is zero comes from comparing with the Poincar\'e duality spin bordism group in which $X$ itself gives a zero bordism.
 Even though this is not really important for our main argument, we next show that $f$ is unique in a certain sense:
In the spin case the normal cobordism group of degree~1 normal maps  is isomorphic to
\[ H_2(X;\Z/2)
\oplus 8\cdot\Z \]
where the $8\cdot\Z$ summand is the difference of ordinary signatures, $\sigma(N)-\sigma(X)$. The first summand vanishes in our  case and so does $\sigma(X)$. By adding copies of Freedman's
$E_8$--manifold to $N$, we may assume that $\sigma(N)=0$. This is just another way of saying that we work with the {\em reduced} normal cobordism group.
 \subsection{The surgery obstruction}\label{sec:surgery}
There is a well defined surgery obstruction $\sigma(f)$ in $L^h_4(\Z[G])$ which in our case is simply the intersection form on $\pi_2N$. The first step is `surgery below the middle dimension' which makes $f$ 2--connected. The surgery obstruction then measures the kernel of $f$ in the middle dimension which in our case is in $\pi_2$.
To make this precise, recall the following definition from  \cite[p.~47]{W99}.
 \begin{definition}
A {\em quadratic form} over $\Z[G]$ is defined to be a triple $(H,\l,\mu)$ with $H$ a free $\Z[G]$--module,
\[ \l\co  H \times H \to \Z[G] \]
a non-singular hermitian form and
\[ \mu\co  H  \to \Z[G]/\langle a-\bar{a} \mid a \in \Z[G]\rangle \]
a quadratic refinement. Here a form is
{\em non-singular} if the induced map $H\to \Hom_{\Z[G]}(H;\Z[G])$ is an isomorphism.
A form isomorphic to a direct sum of
the form $(\Z[G] \cdot e\oplus \Z[G]\cdot f,\l,\mu)$,
where $\l(e,f)=1$ and $\mu(e)=\mu(f)=0$,
is called a {\em hyperbolic form}.
\end{definition}
It should be pointed out that in the oriented case, where the involution on $\Z[G]$ is given by $\bar g:= g^{-1}$, the quadratic refinement $\mu$ is completely determined by the hermitian form $\l$. Its only role is to make sure that $\l$ is {\em even} in the sense that for every $h\in H$ there is an $m \in \Z[G]$ such that
\begin{equation}\tag{even}
\l(h,h) = m + \bar m
\end{equation}
The main examples of hermitian forms come from the intersection form of $4$--manifolds, where
$G=\pi_1W$ and $(H,\l)=(\pi_2W^4,\l_W)$. If $h\in \pi_2W$ is represented by an immersed
$2$--sphere $S$ in $W$, then one can look at the self-intersections $m(S)$ that are related to
the intersection $\l_W$ of $S$ and a parallel push-off by the following formula:
\[
\l_W(S,S) = m(S) + \bar m(S) + e(S)\cdot 1
\]
Here $e(S)\in\Z$ denotes the normal Euler class of $S$.
Note that the coefficient of $1\in G$ in $m(S)$ can be changed arbitrarily by a (nonregular) homotopy, the so called {\em cusp move}. This also changes $e(S)$ so as to make the above formula still hold (the left hand side only depends on the homotopy class of $S$).
It follows that the intersection form $\l$ on $\pi_2$ of a $4$--manifold is even in the above sense if and only if $e(S)$ is an even integer, ie,  the second Stiefel--Whitney class $w_2(S)$ vanishes. That is why we will completely ignore the quadratic refinement $\mu$ in the following, only making sure that our manifolds are spin.
 Note that for general spin $4$--manifolds $W$ with fundamental group $G$, the intersection form $\l_W$ on $\pi_2W$ is not a quadratic form in the above sense. The problems are that in general $\pi_2W$ is not free and  $\l_W$ is not non-singular. However, given a $2$--connected degree~1 map $f\co N\to X$, one can restrict $\l_N$ to the kernel of $f$ on $\pi_2$. Then both of these conditions can be arranged after stabilizing $N$ by copies of $S^2\times S^2$, \cite[p.~26]{W99}. Moreover, if $w_2$ vanishes on $2$--spheres in this kernel then a quadratic refinement exists by the above considerations. If $f$ is a normal map (in addition to having degree~1) then this last condition is obvious.
 In our setting, $\pi_2K(G,1)=0$ so that we automatically work on the kernel of $f$. Moreover, we saw that $N$ is spin so that the intersection form $\l_N$ on $\pi_2N$ is even and hence represents a unique quadratic form, up to the stabilization with hyperbolic forms. This motivates the following definition.
Consider the semigroup of quadratic forms under direct sum.
We say that forms $X_1,X_2$ are equivalent if there exist hyperbolic
forms $H_1,H_2$, such that $X_1\oplus H_1$ and $X_2\oplus H_2$ are
isomorphic. The set of equivalence classes form a group
(cf \cite[p.~249]{R02}), denoted by $L_4^h(\Z[G])$.  As discussed above, $(\pi_2N,\l_N)$ represents an element in this group.
 \begin{theorem} \label{thm:lgroups}
The ordinary signature $\sigma$ induces isomorphisms
\begin{enumerate}
\item
$ L_4^h(\Z) \overset{\cong}{\ra}  \Z, \quad (H,\l,\mu)\mapsto \sigma(\l)/8$
\item $L_4^h(\Z[G])\overset{\cong}{\ra}  L_4^h(\Z)$ for $G=\Z$ or $G=SR$.
\end{enumerate}
\end{theorem}
 The first statement is well-known, see for example \cite[Theorem~13A.1]{W99}. The second statement
is well-known in the case that $G=\Z$ (cf \cite{R98}). We prove the case $G=SR$ in Section
\ref{sec:SR}. In general if $G$ is a ribbon group, then $L_4^h(\Z[G])\to L_4^h(\Z)$ is an isomorphism
if the Whitehead conjecture and the Farrell--Jones conjecture hold. This has been shown in many interesting cases, cf \cite{AFR97} for the
case of knot groups. By the computations in Section~\ref{sec:degree 1}, the surgery obstruction $\sigma(f)=\l_N$ actually only depends on the original data $(M,\varphi)$ if we assume that $\sigma(N)=0$, ie, that it lies in the reduced L--group
\[
\sigma(M,\varphi) = \sigma(f) = \l_N \in \widetilde L^h_4(\Z[G]) := \ker  (L^h_4(\Z[G]) \ra  L^h_4(\Z))
\]
This element is hence an obstruction for finding a h--ribbon for $K$ with epimorphism $\phi\co \pi_1M\onto G$.
If the reduced L--group vanishes then, after further stabilization by $S^2\times S^2$, we may assume that $\pi_2N$ has hyperbolic intersection form.
 \subsection{Constructing the h--ribbon disk} \label{sec:final}
The assumptions of Theorem~\ref{thm:mainG} say that the reduced L--group vanishes and
that the topological surgery sequence $(\mathcal S)$ is exact for our fundamental group. Therefore, one gets a topological manifold $W$ together with a homotopy equivalence $(W,M_K)\to (K(G,1), M_K)$. Hence we found our h--ribbon by Proposition~\ref{prop:topslice} and we finished the proof of Theorem~\ref{thm:mainG}.

\subsection{Proof of Theorem~\ref{thm:main}} \label{sec:main}
 To construct a ribbon disk we use
Theorem~\ref{thm:mainG} that we just finished proving. For $G=SR$, we will show all the required properties in
Section~\ref{sec:SR}. For $G=\Z$ they are easy to check, for example the conditions on $\Hom$
follow from the fact that the Alexander module $H_1(M_K;\Z[\Z])$ is a $\Z[\Z]$--torsion module.
Moreover, \[ \Ext_{\Z[\Z]}^1(H_1(M_K;\Z[\Z]),\Z[\Z]) \cong H_1(M_K;\Z[\Z]) \] so that our (Ext)
condition simply means that the Alexander module (or, equivalently, the Alexander polynomial) of
$K$ is trivial.
 Conversely, assume that there exists a h--ribbon $D$ for $K$ with $\pi_1(N_D)\cong
G$. Then it follows from arguments as in the proof of Theorem 4.4 in \cite{COT03} that the
Blanchfield form $\Bl(G)$ always vanishes on \[ \Ker\{H_1(M_K;\Z[G])\ra H_1(N_D;\Z[G])\}, \] In our
case this group equals $H_1(M_K;\Z[G])$ since $H_1(N_D;\Z[G])=0$. In the Section \ref{sec:BL} we
will recall this Blanchfield form and show in Lemma~\ref{lem:Bl=0} that its vanishing is indeed
equivalent to our (Ext) condition.

\section{Properties of the solvable ribbon group $SR$} \label{sec:SR}
 The purpose of this section
is to prove that the assumptions of Lemma~\ref{lem:pcext} are satisfied for the group $G=SR$:
 \begin{lemma}\label{lem:LSR} The group $G=SR$ satisfies $H_3(G)=0$ and $H^i(G;\Z[G])=0$ for $i>2$. Furthermore, for any knot $K$ and any epimorphism $\pi_1(M_K) \onto G$ one has
\[
 \Hom_{\Z[G]}(H_i(M_K;\Z[G]),\Z[G])=0\quad  \text{ for }\quad  i=1,2.
\]
\end{lemma}
A group $G$ is called poly-torsion-free-abelian (PTFA) if there exists a
filtration \[ 1=G_0 \vartriangleleft G_1\vartriangleleft \dots \vartriangleleft
G_{n-1}\vartriangleleft G_n=G \] such that $G_{i}/G_{i-1}$ is torsion free abelian. As we pointed
out in Section~\ref{sec:ribboncomp} the only solvable ribbon groups are $\Z$ and $SR$ and clearly
both are PTFA. If $G$ is a PTFA group then the group ring
$\Z[G]$ satisfies the {\em Ore condition} and therefore has an {\em Ore localization}, cf \cite[Proposition~2.5]{COT03}
\[
Q(G):=\Z[G](\Z[G]\sm 0)^{-1}
\]
which is a skew field. Therefore, every $Q(G)$--module is free and its rank is well-defined.
As a localization, $Q(G)$ is flat as a right (and left) $\Z[G]$--module. Moreover, the $\Z[G]$--torsion part of a module is a submodule (unlike for general non-commutative rings).
 Let $K$ be a knot and $\pi_1(M_K)\onto G$ be an epimorphism onto a PTFA group. From
\cite[Proposition~2.11]{COT03} it follows that
\[
H_1(M_K;\Z[G]) \otimes_{\Z[G]} Q(G)\cong H_1(M_K;Q(G))=0,
\]
ie, $H_1(M_K;\Z[G])$ is $\Z[G]$--torsion. In particular, the required condition
\[ \Hom_{\Z[G]}(H_1(M_K;\Z[G]),\Z[G])=0 \] in Lemma~\ref{lem:LSR} follows. Furthermore,
\[
H_0(M_K;Q(G))\cong H_0(M_K;\Z[G])\otimes_{\Z[G]}Q(G)\cong\Z \otimes_{\Z[G]}Q(G)=0 \]
and similarly $H_3(M_K;Q(G))=0$ since $H_3(M_K;\Z[G])=0$.
The Euler characteristics of $M_K$ with $\Z$ and $Q(G)$ coefficients agree by the usual argument
and hence\[ \chi_{Q(G)}(H_*(M_K;Q(G)))=\chi_\Z(H_*(M_K;\Z))=0. \] It follows that \[
H_2(M_K;\Z[G])\otimes_{\Z[G]}Q(G)\cong H_2(M_K;Q(G))=0, \]
Therefore $H_1(M_K;Q(G))=0$ and hence $\Hom_{\Z[G]}(H_2(M_K;\Z[G]),\Z[G])=0$ because $H_2(M_K;\Z[G])$ is again
$\Z[G]$--torsion. So another required condition in Lemma~\ref{lem:LSR} follows: \[
\Hom_{\Z[G]}(H_2(M_K;\Z[G]),\Z[G])=0. \]
 The following result finishes up our task of showing $SR$ has all the required properties.
\begin{lemma} \label{lem:LSR2} The group $G=SR$ is aspherical, in particular
$\Ext_{\Z[G]}^i(\Z,\Z[G])=0$ for $i>2$ and $H_3(G)=0$.
\end{lemma}
 \begin{proof}
Let $D$ be a ribbon disk with  complement  $N_D$ and fundamental group $\pi_1(N_D)=G=SR$. Recall
that by Alexander duality $H_2(N_D)=0$. We will show that  $N_D$ is a $K(G,1)$. First note that
$N_D$ is homotopy equivalent to a 2--complex $N$. By the Hurewicz theorem it suffices to prove the
vanishing of $\pi_2(N) \cong H_2(\ti{N})$.
 Let $Q(G)$ be the Ore localization of $\Z[G]$.
By the same arguments as above for $M_K$, involving Euler characteristics, one shows that
\[
H_2(\ti{N})\otimes_{\Z[G]}Q(G)\cong H_2(N;\Z[G])\otimes_{\Z[G]} Q(G)=0.
\]
Now consider the following commutative diagram
\[ \ba{rcccc}
0&\rightarrow &H_2(\ti{N})&\hookrightarrow& C_2(\ti{N})\\
&&\downarrow&&\downarrow \\
0&\to &H_2(\ti{N})\otimes_{\Z[G]} Q(G)&\hookrightarrow& C_2(\ti{N})\otimes_{\Z[G]} Q(G).\ea \]
The second vertical map is injective since $C_2(\ti{N})$ is a free $\Z[G]$--module, hence the
first vertical map is injective as well, hence $H_2(\ti{N})=0$.
\end{proof}
 \begin{remark}
We could have used the well known fact that a one-relator group is aspherical if the
relation is not a proper power (so in particular the group is torsion free).
However, a similar argument as above will also be used in the proof of Lemma \ref{lem:hig}.
Locally indicable ribbon groups are also known to be aspherical.
\end{remark}
 The following result is not needed in the rest of the paper but we include it for the interested
reader. The nontriviality of the homology group in question made it clear that some more naive
formulations of the  Poincar\'e duality condition must fail. Recall from
Section~\ref{sec:ribboncomp} that $SR$ has the presentation \[ SR=\langle a,c \mid
aca^{-1}=c^2\rangle. \]
 \begin{lemma} \label{lem:hig}
Let $G:=SR$, then $H^2(G;\Z[G])$ maps onto $\Z[1/2]$. More precisely, as a right $\Z[G]$--module, $H^2(G;\Z[G])$ is the quotient of $\Z[G]$ by the right ideal generated by
\[
1- aca^{-1} \quad \text{ and }\quad a-aca^{-1}c^{-1}-aca^{-1}c^{-2}.
\]
Dividing $\Z[G]$ by the {\em two-sided} ideal generated by the same elements gives the ring $\Z[1/2]$.
\end{lemma}
 \begin{proof} The proof of Lemma~\ref{lem:LSR} shows  that the 2--complex associated to the above presentation is a $K(G,1)$. Therefore, we can calculate $H^2(G;\Z[G])$ from the short resolution
\[
0\ra \Z[G] \ra \Z[G]^2 \ra \Z[G] \overset{\epsilon}{\ra} \Z \ra 0
\]
where the first map is given
by the Fox derivatives $\partial_a$ and $\partial_c$ of the relation $aca^{-1}c^{-2}$. A straightforward
calculation using Fox derivatives `from the left' gives
\[ \partial_a = 1- aca^{-1} \quad \text{ and } \quad
\partial_c=a-aca^{-1}c^{-1}-aca^{-1}c^{-2}
 \]
If we now decide to quotient by the two-sided ideal generated by the first
relation, we find that the defining relation forces $c=1$ and hence we are left with $\Z[\Z]$,
generated by $a$. The second relation then introduces the relation $a=2$ into this commutative ring
which gives the  ring $\Z[1/2]$.
 \end{proof}

\subsection*{Computation of $L_4^h(\Z[SR])$}

\begin{lemma} \label{lem:LSR4}
The  inclusion map $\Z \to \Z[\Z\ltimes \Z[1/2]]$ induces an
isomorphism $L_4^h(\Z) \to L_4^h(\Z[\Z\ltimes \Z[1/2]])$.
\end{lemma}
 In order to prove this lemma, we need a theorem by Ranicki,
which in turn needs the notion of $L_n$--groups for any $n$. We refer to \cite{R92}
for the definition of these groups.
 We recall that the Whitehead group $\Wh(G)$ for a group $G$
is defined as $\mbox{Wh}(G)=K_1(\Z[G])/\pm G$
(cf \cite[p.~172]{R02} for a definition of $K_1(\Z[G])$).
\begin{theorem}\label{thmranicki}{\rm\cite{R73}}\qua
Let $G$ be a group with $\Wh(G)=0$ and let $\a\co \Z [G]\to \Z [G]$
be an automorphism, then there exists
a long exact sequence
\[ \dots \to L_n(\Z[G])\xrightarrow{1-\a_*} L_n(\Z[G])\to
L_n(\Z[G]^{\a}[t,t^{-1}])\to  L_{n-1}(\Z[G])
\to \dots, \]
where $\Z[G]^{\a}[t,t^{-1}]$ denotes the twisted Laurent ring,
ie, for $a_1,a_2\in \Z[G]$ we have
$a_1t^{n_1}\cdot a_2t^{n_2}=a_1\a(a_2)^{n_1}t^{n_1+n_2}$.
\end{theorem}
 \begin{remark}
Since $Wh(G)=0$, we don't have to distinguish between $L^s$ and  $L^h$. We will therefore henceforth drop any decorations.
\end{remark}

We'll make  use of the fact that $\Z[1/2]$ is the direct limit
of the direct system $\Z\xrightarrow{\cdot 2}\Z
\xrightarrow{\cdot 2}\dots $. We'll therefore write
$\Z[1/2]=\underset{\to}{\lim} \, \Z$.
Denote the map $\cdot 2\co \Z[1/2]\to \Z[1/2]$ by $\a$; then $\Z[\Z\ltimes
\Z[1/2]]\cong \Z[\Z[1/2]]^{\a}[t,t^{-1}]$.
 \begin{claim}
$\quad\mbox{Wh}(\Z[1/2])=0.$
\end{claim}
 The determinant map $\det\co  K_1(\Z[\Z])\to \Z \times \{\pm 1\}$ is an
isomorphism
(cf \cite[p.~172]{R02}). The $K_1$ functor commutes with direct limits and
$\a$ commutes with the determinant maps, hence
$\det\co  K_1(\Z[\Z[1/2]])\to \Z[\Z[1/2]]\times \{\pm 1\}$ is an isomorphism
as well.
It now follows immediately that $ \mbox{Wh}(\Z[1/2])=0$.

By Theorem \ref{thmranicki} there therefore exists an exact sequence
\[ \ba{rcl} \dots \to &  L_n(\Z[\Z[1/2]])&\xrightarrow{1-\a_*} L_n(\Z[\Z[1/2]])\to
L_n(\Z[\Z[1/2]]^{\a}[t,t^{-1}])\to \\ \to & L_{n-1}(\Z[\Z[1/2]])&
\xrightarrow{1-\a_*}  \dots \ea \]
 It is clear that the proposition follows once we prove the following claim.
\begin{claim}
\bn
\item $L_4(\Z)\to L_4(\Z[\Z[1/2]]$ is an isomorphism,
\item $\a_*\co L_4(\Z[\Z[1/2]])\to L_4(\Z[\Z[1/2]])$ is the identity map,
\item $L_3(\Z[\Z[1/2]])=0.$
\en
\end{claim}
 We recall that (cf \cite{W99})
\[ L_n(\Z)=\left\{ \ba{rcl} 0 &\mbox{ if}& n \equiv 1 \mod 2 \\
\Z &\mbox{ if}& n \equiv 0 \mod 4 \\
\Z/2 &\mbox{ if}& n \equiv 2 \mod 4. \ea \right. \] From  Theorem \ref{thmranicki} it follows
that $L_4(\Z])\xrightarrow{\cong} L_4(\Z[\Z])$ via the map induced by the inclusion $\Z\to
\Z[\Z]$. The inverse $L_{4}(\Z[\Z])\to L_{4}(\Z)$ is given by tensoring with $\Z$, considered as
a $\Z[\Z]$--module via the map $t\mapsto 1$. In particular the map $\Z \xrightarrow{\cdot 2}\Z$
induces the identity on the $L_4$ group of the group ring, ie, on $L_4(\Z[\Z])$. From the fact
that the $L$--functor commutes with direct limits, we immediately get the first two statements.
 From \cite[p.~181]{W99} we get the following commutative diagram
\[ \ba{ccc} L_3(\Z[\Z])&\to &L_3(\Z[\Z/2]) \\ \hbox{\scriptsize $\a_*$} \downarrow&&
\downarrow \hbox{\scriptsize $\a_*$} \\ L_3(\Z[\Z])&\to &L_3(\Z[\Z/2]),\ea \] where $\a$ denotes multiplication by
two. In particular $\a_*\co L_3(\Z[\Z/2])\to$\break $L_3(\Z[\Z/2])$ factors through $L_3(\Z)=0$. Hence
$\a_*\co  L_3(\Z[\Z])\to  L_3(\Z[\Z])$ is the zero map, taking direct limits we see that
$L_3(\Z[\Z[1/2]])=0$. This completes the proof of Lemma~\ref{lem:LSR4}.

\section{Non-commutative Blanchfield forms} \label{sec:BL}
Let $G$ be a PTFA--group and $\varphi\co M_K\onto G$ a homomorphism. Note that the involution on
$\Z[G]$ extends to an involution on $Q(G)$ and on $Q(G)/\Z[G]$. Let $\pi:=\pi_1(M_K)$. Note that
the map
\[ \ba{rcl} \hspace{-5pt}\Hom_{\Z[\pi]}(C_*(\ti{M}_K),Q(G)/\Z[G]) &\to&
\Hom_{\Z[G]}(C_*(\ti{M}_K)\otimes_{\Z[\pi]} \Z[G],Q(G)/\Z[G])\hspace{-3pt}\\
\phi&\mapsto& (c\otimes f\mapsto \overline{\phi(c)}f )\ea \] is well--defined and induces an
evaluation map
\[ H^1(M_K;Q(G)/\Z[G])
\ra \Hom_{\Z[G]}(H_1(M_K;\Z[G]),Q(G)/\Z[G]).  \] Now consider
\[ \ba{rcl} H_1(M_K;\Z[G])&\overset{\cong}{\ra}& H^2(M_K;\Z[G])\overset{\cong}{\longleftarrow} H^1(M_K;Q(G)/\Z[G])\\
&\ra& \Hom_{\Z[G]}(H_1(M_K;\Z[G]),Q(G)/\Z[G])\ea
\]
ie, the composition of Poincar\'e duality, the inverse of the Bockstein homomorphism for the
coefficient sequence $0\to \Z[G]\to Q(G)\to Q(G)/\Z[G]\to 0$ and the above evaluation map.

Note that the second map is an  isomorphism since the homology and hence the cohomology with $Q(G)$
coefficients vanishes.
 This map defines the hermitian {\em Blanchfield pairing},  cf \cite[Theorem~5.1]{D86}
\[
\Bl(G)\co  H_1(M_K;\Z[G])\times H_1(M_K;\Z[G])\ra Q(G)/\Z[G].
\]

\begin{lemma} \label{lem:Bl=0}
For $G=\Z$ or $SR$, the condition {\em (Ext)} from Theorem~\ref{thm:main} is equivalent to the vanishing of the Blanchfield form $\Bl(G)$.
\end{lemma}
 The remark after Lemma~\ref{lem:pc} shows that the vanishing of the Blanchfield pairing $\Bl(G)$ is a weaker statement than
the vanishing of the corresponding homology group $H_1(M_K;\Z[G])\cong H^2(M_K;\Z[G])$. In
particular, this shows that $\Bl(G)$ will in general be singular.

\begin{proof}[Proof of Lemma~\ref{lem:Bl=0}]
From the long exact $\Ext$--sequence corresponding to
\[
0\ra \Z[G]\ra Q(G)\ra Q(G)/\Z[G]\ra 0
\]
it follows that
\[ \ba{rcl} \Ext_{\Z[G]}^1(H_1(M_K;\Z[G]),\Z[G])&\cong &\Ext_{\Z[G]}^0(H_1(M_K;\Z[G]),Q(G)/\Z[G])\\
& \cong &\Hom_{\Z[G]}(H_1(M_K;\Z[G]),Q(G)/\Z[G])\ea \] since
$\Ext_{\Z[G]}^i(H_1(M_K;\Z[G]),Q(G))=0$ for all $i$ (cf \cite[Remark~2.8]{COT03}). We are now
done once we show that the homomorphism
\begin{equation}  \label{equ:surj} H^1(M_K;Q(G)/\Z[G])\ra \Hom_{\Z[G]}(H_1(M_K;\Z[G]),Q(G)/\Z[G])\end{equation} in the
definition of $\Bl(G)$ is surjective. As in the proof of Lemma \ref{lem:pcext} we have an
isomorphism \[ H^1(M_K;Q(G)/\Z[G])\cong
H_1(\Hom_{\Z[G]}(C_*(\ti{M}_K)\otimes_{\Z[\pi_1(M)]}\Z[G],Q(G)/\Z[G])).\] Now applying the UCSS
  we see
that the evaluation map (\ref{equ:surj}) is surjective if $\Ext_{\Z[G]}^2(\Z,Q(G)/\Z[G])=0$.
 Let $X=K(SR,1)$. By Lemma \ref{lem:LSR2} we can assume that $X$ is a 2--complex. Then
\[ \Ext_{\Z[G]}^2(\Z,Q(G)/\Z[G])\cong
H^2(G;Q(G)/\Z[G]) \cong H^2(X;Q(G)/\Z[G]). \] Note that $H^i(X;Q(G))\cong H_i(X;Q(G))\cong
H_i(X;\Z[G])\otimes_{\Z[G]} Q(G)=0$ for $i=2,3$, in particular it follows that
\[
H^2(X;Q(G)/\Z[G])=H^3(X;\Z[G])=0
\]
because $X$ is a 2--complex.
\end{proof}
 Since $\Bl(\Z)$ is nonsingular for any knot $K$, it vanishes if and only if the Alexander module
itself vanishes. This in turn is equivalent to $K$ having trivial Alexander polynomial and to our
condition (Ext) for  $G=\Z$.

\section{The satellite construction} \label{sec:satellite} 
\begin{figure}[ht!]\small\anchor{figsat} 
\begin{center}
\psfrag {C}{$C$}
\psfrag {K}{$K$}
\psfrag {A}{$A$}
\psfrag {S}{$S$}
\includegraphics[scale=0.31]{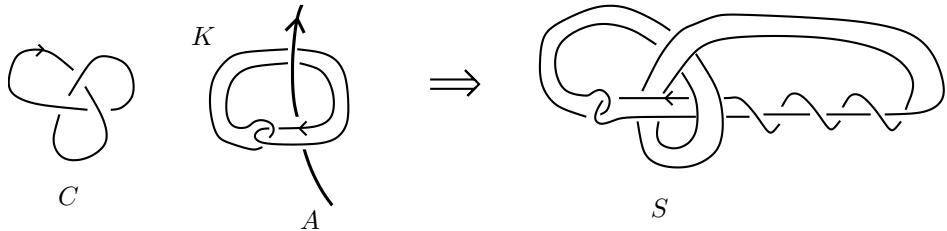}
\end{center}
\caption{Satellite construction with $C$ the trefoil, $K$ the unknot}
\label{figsat}
\end{figure}

Let $K,C$ be knots. Let $A\subset S^3\sm
K$ be a curve, unknotted in $S^3$. Then $S^3 \smallsetminus \nu A$ is a solid torus. Let
$\psi\co \partial(\overline{\nu A})\to \partial(\overline{\nu C})$ be a diffeomorphism which sends a meridian of $A$ to a
longitude of $C$, and a longitude of $A$ to a meridian of $C$. The space
\[ (S^3\sm \nu A)\cup_{\psi} (S^3\sm \nu C) \] is a 3--sphere and the image of $K$ is denoted by
$S=S(K,C,A)$. We say $S$ is the satellite knot with companion $C$, orbit $K$ and axis $A$. Note
that we replaced a tubular neighborhood of $C$ by a knot in a solid torus, namely $K\subset
S^3\sm \nu A$. \figref{figsat} shows that taking the Whitehead double of a knot is an
example for the satellite construction.
 It is easy to see that \[ M_S=(M_K\sm \nu A) \cup_{\partial(\overline{\nu A})}  (S^3\sm \nu C).
\] Obstruction theory shows that the map $\psi^{-1}\co \partial(\overline{\nu C})\to \partial(\overline{\nu A})$ can be
extended to a map $f\co S^3\sm \nu C\to \overline{\nu A}$. Combining with the injection
$M_K\sm \nu A \to M_K$ this defines a map $M_S\to M_K$. Let $\varphi\co \pi_1(M_K)\to G$ be
a homomorphism onto a torsion-free group. Then we get an induced map $\pi_1(M_S)\to G$.
 \begin{lemma} \label{lemmahofsat} If $\varphi(A)=e$, then $H_1(M_S;\Z[G])\cong H_1(M_K;\Z[G])$,
otherwise
\[ H_1(M_S;\Z[G])\cong H_1(M_K;\Z[G])\oplus H_1(M_C;\Z[\Z])\otimes_{\Z[\Z]} \Z[G],\]
where $\Z[G]$ is a $\Z[\Z]\cong \Z[t,t^{-1}]$ module via $t\mapsto \varphi(A)$.
\end{lemma}
 \begin{proof}
Consider the following commutative diagram of Meyer--Vietoris exact sequences (with
$\Z[G]$--coefficients understood)
\[ \ba{ccccccccc}
\hspace{-0.05cm}H_1(\partial(\overline{\nu A}))\hspace{-0.05cm} & \hspace{-0.0cm}\longrightarrow
\hspace{-0.0cm}&\hspace{-0.05cm}H_1(M_K\sm \nu A)\hspace{-0.05cm} &\hspace{-0.0cm}\oplus
\hspace{-0.0cm} & \hspace{-0.05cm}H_1(S^3\sm \nu C)\hspace{-0.05cm}&
\hspace{-0.0cm}\longrightarrow \hspace{-0.0cm}& \hspace{-0.05cm}H_1(M_S)\hspace{-0.05cm}&\hspace{-0.0cm}\longrightarrow\hspace{-0.0cm}&\hspace{-0.05cm} 0 \\
\downarrow&&\downarrow&&\downarrow\! \raise 3pt\hbox{\scriptsize $f$}&&\downarrow&\\
\hspace{-0.05cm}H_1(\partial(\overline{\nu A}))\hspace{-0.05cm}
&\hspace{-0.0cm}\longrightarrow\hspace{-0.0cm} &\hspace{-0.05cm}H_1(M_K\sm \nu
A)\hspace{-0.05cm} &\hspace{-0.0cm}\oplus\hspace{-0.0cm}&\hspace{-0.05cm} H_1(\overline{\nu
A})\hspace{-0.05cm} &\hspace{-0.0cm}\longrightarrow \hspace{-0.0cm}&
\hspace{-0.05cm}H_1(M_K)\hspace{-0.05cm}&\hspace{-0.0cm}\longrightarrow\hspace{-0.0cm}& 0\ea
\]
If $\varphi(A)=e$, then
the coefficient systems for $\overline{\nu A}$ and for $S^3\sm \nu C$
are trivial, hence
\[
\ba{rcl} H_1(\overline{\nu A};\Z[G])&\cong &H_1(\overline{\nu A};\Z)\otimes_{\Z} \Z[G]\\
&\cong &\Z\otimes_{\Z}\Z[G]\\
&\cong & H_1(S^3\sm \nu C;\Z)\otimes_{\Z}\Z[G] \\
&\cong &H_1(S^3\sm \nu C;\Z[G]). \ea \]
This immediately implies that $H_1(M_S;\Z[G])\cong
H_1(M_K;\Z[G])$.
 If $\varphi(A)\ne e$, then  $H_1(\overline{\nu A};\Z[G])=0$ since $\varphi$ is an element of
infinite order since $G$ is torsion free. Furthermore $H_1(\partial(\overline{\nu A}))$ is a free
$\Z[G]$--module on the meridian of $A$, which gets mapped to the longitude in $S^3\sm \nu C$ which
is zero in $H_1(S^3\sm \nu C;\Z[\Z])$. From the above commutative diagram it now follows that \[
H_1(M_S;\Z[G])\cong H_1(M_K;\Z[G])\oplus H_1(S^3\sm \nu C;\Z[\Z])\otimes_{\Z[\Z]} \Z[G]\] since
$H_1(S^3\sm \nu C;\Z[G])\cong H_1(S^3\sm \nu C;\Z[\Z])\otimes_{\Z[\Z]}\Z[G]$. The lemma now follows
from $H_1(S^3\sm \nu C;\Z[\Z])\cong H_1(M_C;\Z[\Z])$. \end{proof}

\section{Examples: Satellite knots of $6_{1}$} \label{sec:61} 
\begin{figure}[ht!]\small\anchor{exknots} 
\begin{center} 
\psfrag {K1}{$K_1$}
\psfrag {K2}{$K_2$}
\psfrag {K3}{$K_3$}
\includegraphics[scale=0.25]{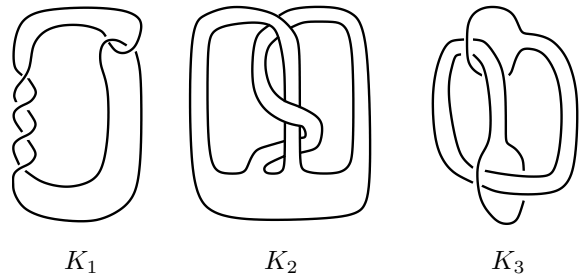}
\caption{Three projections of the knot $6_1$} \label{exknots} 
\end{center}
 \end{figure}

\begin{figure}[ht!]\small\anchor{ribbonfigure}
\begin{center}\includegraphics[scale=0.25]{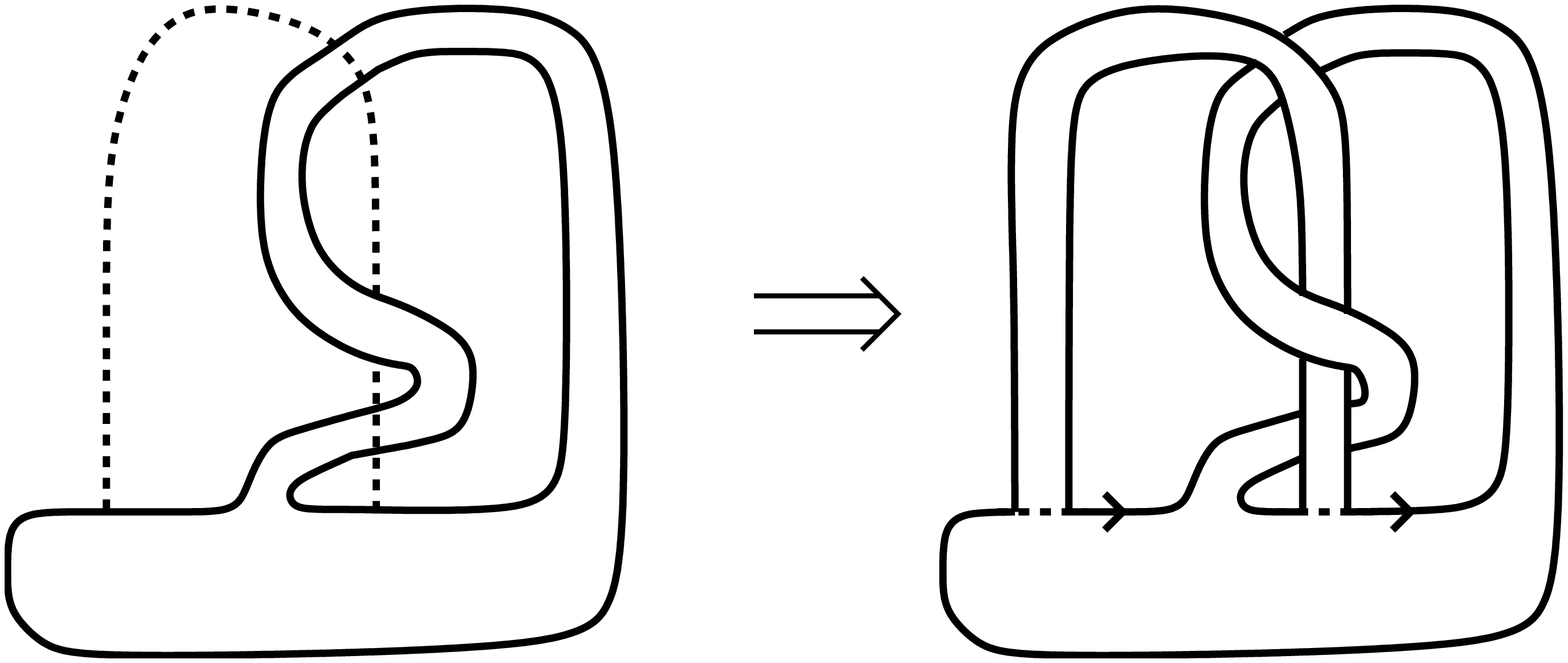} 
\caption{A ribbon disk for $6_1$}
\label{ribbonfigure} \end{center} 
\end{figure}

\begin{figure}[ht!]\small\anchor{knot61} 
\begin{center}
\psfrag {a}{$a$}
\psfrag {b}{$b$}
\psfraga <0pt,-15pt> {al}{$\alpha$}
\psfraga <2pt,-18pt> {be}{$\beta$}
\includegraphics[scale=0.35]{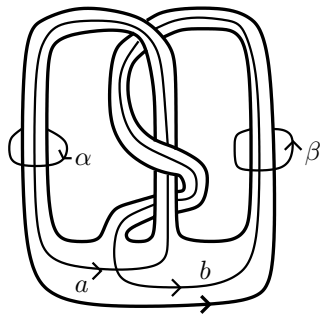} 
\caption{Seifert surface for $6_1$} \label{knot61}
\end{center} \end{figure}

In \figref{exknots} we see three
projections of the knot $6_1$, the first one having the minimal crossing number~6.
Indeed, the isotopy between $K_1$ and $K_3$ is shown in \cite{CG86}. The isotopy between $K_2$ and $K_3$
follows from Figures~\fref{fig:61-C} and \fref{ribbonfigure}, which shows that both knots are formed
by band connected sum of two trivial knots along isotopic bands. Note that the two trivial circles
in \figref{ribbonfigure} bound disjoint disks in $S^2$. 
 Now consider the knot $6_1$ with the Seifert surface $F$ in \figref{knot61}. With the given
basis $a,b$ for $H_1(F)$ we get \[ A(6_1)=\bp 0&2\\  1&0 \ep \] which then shows that
$\Delta_{6_1}(t)=\det(At-A^t)=(t-2)(t^{-1}-2)$. 

\subsection*{Examples of h--ribbon knots} We first prove a general result. Let $R$ be a ribbon knot
with ribbon disk $D_R$, such that $\pi_1(D^4\sm D_R)\cong SR$. Denote the induced map
$\pi_1(M_K)\onto G:=SR$ by $\varphi$. Let $C$ be any knot and $A\subset S^3\sm K$ the unknot in $S^3$.
 \begin{proposition} \label{propsatribbon}
If $\varphi(A)=e$ then the satellite knot $S=S(R,C,A)$ is h--ribbon.
\end{proposition}
 \begin{proof} By Lemma~\ref{lemmahofsat} there exists a map $\pi_1(M_S)\onto
G$ with $H_*(M_S;\Z[G])\cong H_*(M_R;\Z[G])$. From Theorem~\ref{thm:main} it follows that
$\Ext^1_{\Z[G]}(H_1(M_R;\Z[G]),\Z[G])$ $=0$ and hence $(M_S,\varphi)$ also satisfies
$\Ext^1_{\Z[G]}(H_1(M_S;\Z[G]),\Z[G])=0$. The proposition is implied by
Theorem~\ref{thm:main}.
\end{proof}
This result is similar to the well-known fact that the Whitehead double
of any knot (which is the satellite of the unknot)
is topologically slice. This is an immediate corollary from Freedman's slicing theorem
and Lemma~\ref{lemmahofsat}.

\begin{figure}[ht!]\small\anchor{knot61map} 
\begin{center} 
\psfrag {a}{$a$}
\psfrag {b}{$b$}
\includegraphics[scale=0.3]{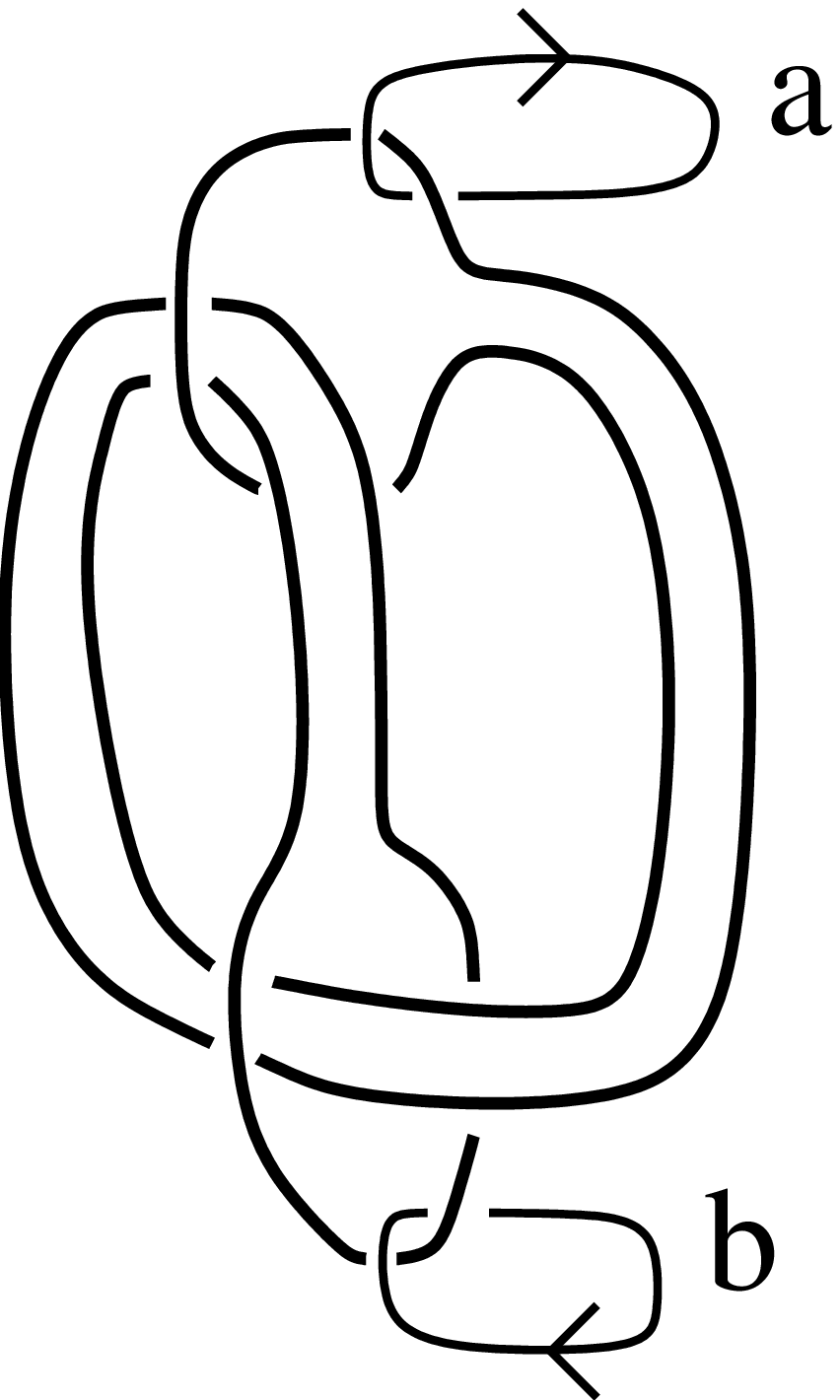}
\caption{The map $\pi_1(S^3\sm K)\to \langle a,b\mid a^{-1}ba^{-1}bab^{-1} \rangle$}
\label{knot61map} \end{center} \end{figure}

 Now turn back to the study of $K=6_1$. From the discussion in Section~\ref{sec:ribboncomp} we know
that $K$ has a ribbon disk $D$ such that $\pi_1(N_D)\cong \langle a,b\mid a^{-1}ba^{-1}bab^{-1}
\rangle =SR$. The group $\pi_1(M_K)\cong \pi_1(S^3\sm K)/\langle \mbox{longitude}\rangle$ is
generated by the meridians of $K_3$. From \figref{figkirby} it follows that the map
\[
\varphi\co \pi_1(M_K)\to \pi_1(N_D)\cong \langle a,b\mid a^{-1}ba^{-1}bab^{-1} \rangle=G \]
is given by
the map indicated in \figref{knot61map}.
The image of the other meridians is determined by the
image of the two given meridians.

\begin{figure}[ht!]\small\anchor{examplea} 
\begin{center} 
\psfrag {A}{$A$}
\psfrag {S}{$S$}
\begin{tabular}{ccc} \includegraphics[scale=0.2]{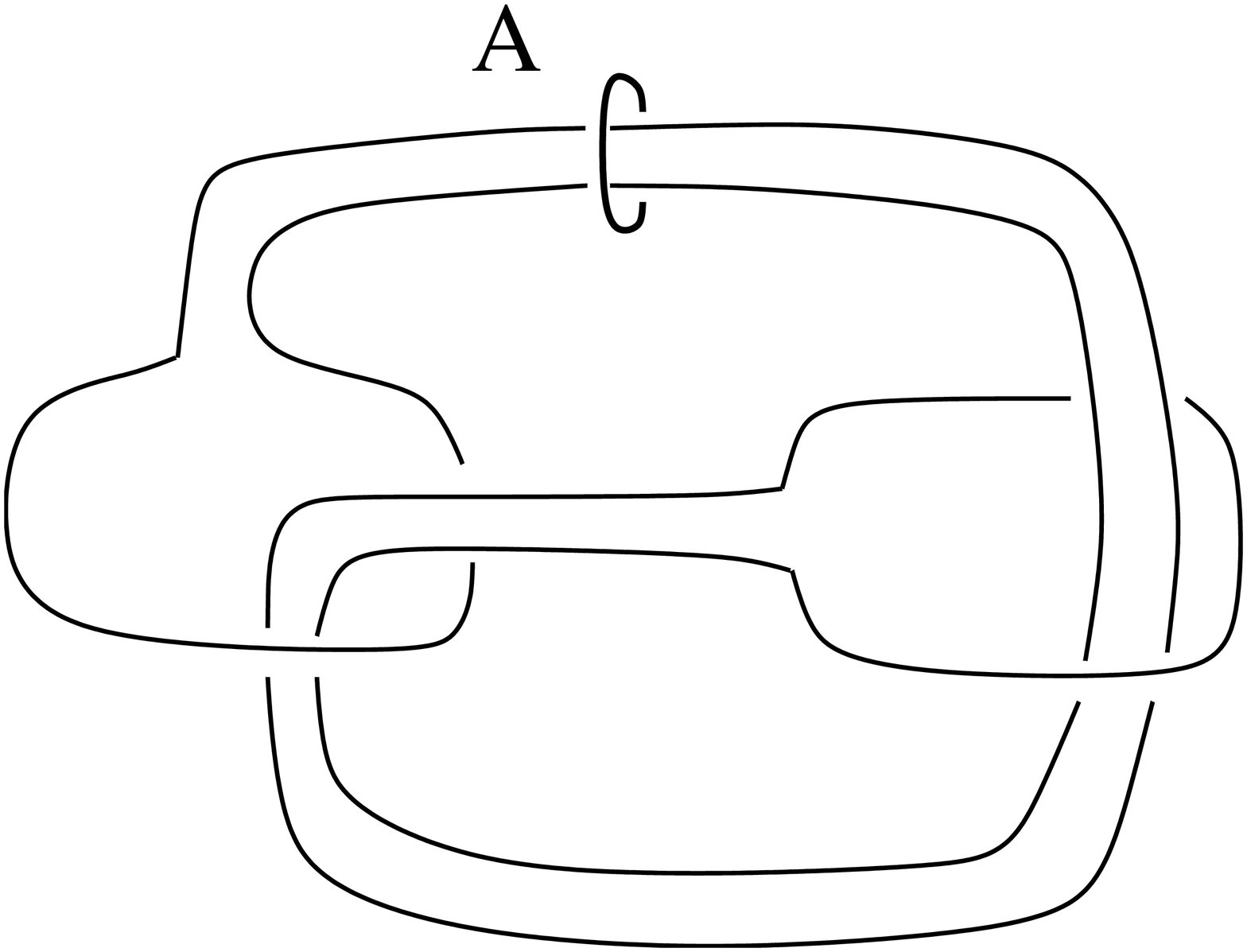}&$$&
\includegraphics[scale=0.07]{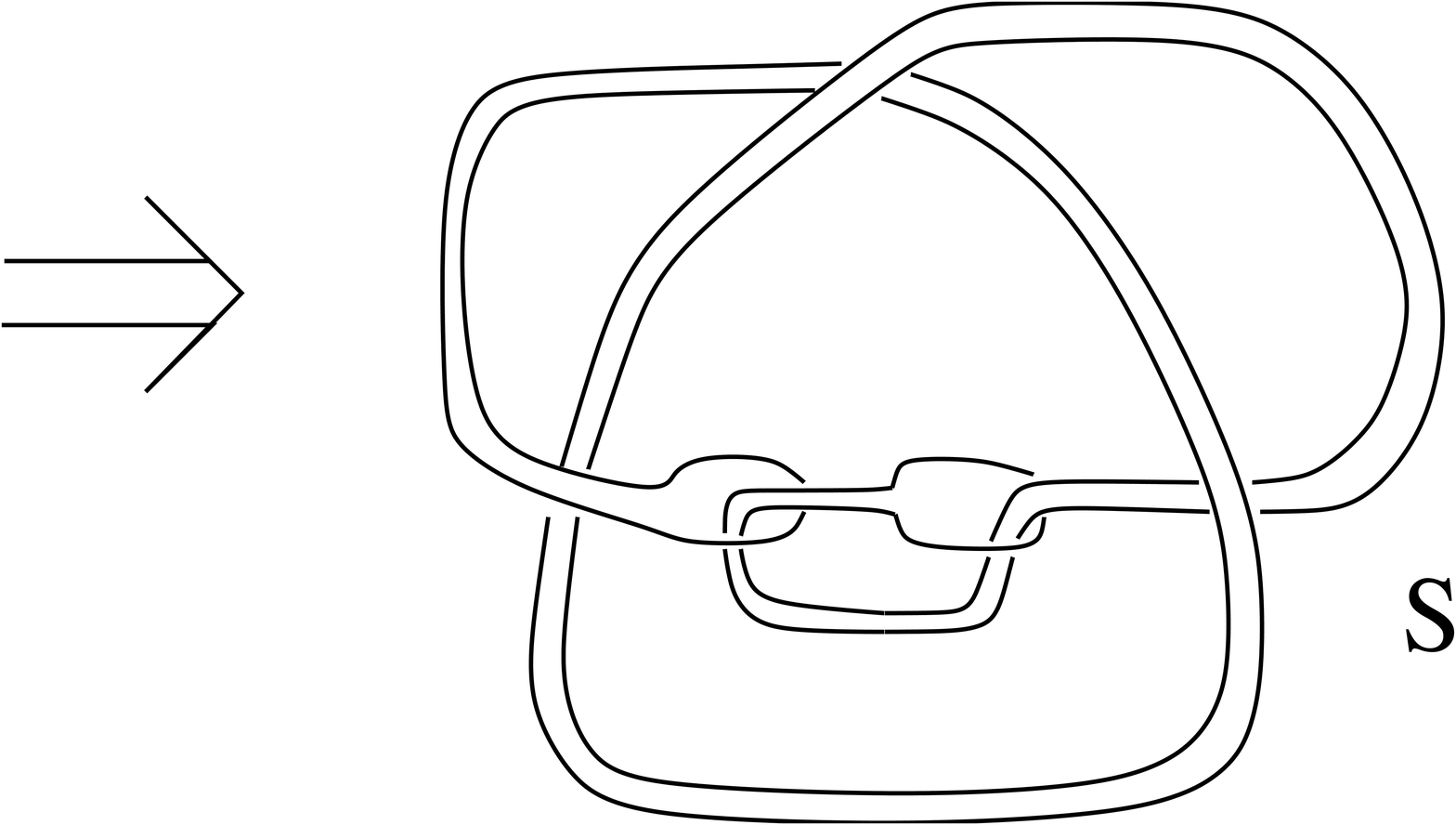}
\end{tabular}
\caption{Ribbon knot with choice of $A$ and the satellite knot
$S(6_1,\mbox{trefoil},A)$}
\label{examplea} \end{center} \end{figure}

As an example take $K=6_1$, $A$ as in \figref{examplea} and $C$ the trefoil knot.
We have $\varphi(A)=e\in G$
as one can easily verify. Therefore $S$ is h--ribbon by Proposition~\ref{propsatribbon}. But in fact
$S$ is trivially ribbon already, since the ribbon disk of $6_1$ `survives' the satellite
construction. We therefore need a more subtle choice of $A$. Now consider \figref{example2}.
One can easily verify that $\varphi(A)=e$. The
knot $6_1$ as a knot in the torus $S^3\sm \nu A$ is given in \figref{example2} on the
right. Wrapping this knot around say the trefoil gives a h--ribbon knot, which we conjecture not to
be smoothly ribbon.
 \begin{remark} It follows immediately from arguments as in Lemma~\ref{lemmahofsat} that
\[
 (H_1(M_K;\Z[\Z]),\Bl(\Z))\cong (H_1(M_S;\Z[\Z]),\Bl(\Z)).
 \]
  In particular $K$ and $S$ have the same
abelian invariants. Livingston   \cite{L02} showed that there exists a knot  $\ti{K}$ with
\[ (H_1(M_K;\Z[\Z]),\Bl(\Z))\cong (H_1(M_{\ti{K}};\Z[\Z]),\Bl(\Z))\] which is not topologically
slice (cf also \cite{K05}). This shows that the vanishing of $\Bl(SR)$, which is a condition
on the non-commutative homology of $K$, can not be reduced to a condition on abelian
invariants.
\end{remark}

\subsection*{Non-existence of certain ribbon disks}

Now consider the knot $6_1$ with the
Seifert surface $F$ in \figref{knot61}.
For given knots  $C_{\a},C_{\b}$ consider $S(K,C_{\a},C_{\b},\a,\b)$, ie, the result
of applying the satellite construction twice. This means that we tie knots $C_{\a}$ and $C_{\b}$ into the bands $\a$ and $\b$.
 \begin{proposition}
If $\Delta_{C_{\a}}(t)\ne 1, \Delta_{C_{\b}}(t)\ne 1$,
then $S$ has no h--ribbon  with
fundamental group $SR$.
\end{proposition}
 \begin{proof} Assume that $S$ has in fact a h--ribbon $D$ with fundamental group
$G:=SR=\Z \ltimes \L/(t-2)$. Then $\Ker\{H_1(M_S;\Z[\Z])\to H_1(N_D;\Z[\Z])\}$ is a metabolizer for
$\Bl(\Z)$ (cf \cite{F04}).
 In the following we will write $\L:=\Z[t,t^{-1}]\cong \Z[\Z]$. Note that $\a,\b$ in
\figref{knot61} lift to elements $\ti{\a}, \ti{\b}$ in $H_1(M_{K};\L)$, in fact \[
H_1(M_{K};\L)\cong (\L\ti{\a} \oplus \L\ti{\b})/(At-A^t)\cong \L\ti{\a}/(2t-1)\oplus
\L\ti{\b}/(t-2).\] Furthermore the Blanchfield pairing $\Bl(\Z)$ is given by the matrix
$(t-1)(At-A^t)^{-1}$. It is easy to see that the submodule generated by $\ti{\a}$ respectively by
$\ti{\b}$ are the only two metabolizers for $\Bl(\Z)$. In particular the map  $\pi\co =\pi_1(M_S)\to
\pi_1(N_D)$ is up to automorphism of $G$ either of the form \[ \varphi_{\a}\co  \pi_1(M_S)\to
\pi/\pi^{(2)} \cong \Z \ltimes H_1(M_S;\L) \to \Z \ltimes \L\ti{\a}/(2t-1) \to SR \] or it is of
the same form with $\ti{\a}$ replaced by $\ti{\b}$. We denote this homomorphism by $\varphi_{\b}$.
By Theorem~\ref{thm:main} we get $\Ext_{\Z[G]}^1(H_1(M_S;\Z[G]),\Z[G])=0$ with $G$--coefficients
induced by $\varphi_{\a}$ or by $\varphi_{\b}$.
 Now consider coefficients induced by $\varphi_{\a}$. From Lemma~\ref{lemmahofsat} it follows that
\[ H_1(M_S;\Z[G])\cong H_1(M_K;\Z[G])\oplus H_1(M_{C_{\a}};\Z[\Z])\otimes_{\Z[\Z]}\Z[G].\] We
compute
\[ \ba{cl} &\Ext^1_{\Z[G]}(H_1(M_S;\Z[G]),\Z[G])\\
\cong &\Ext^1_{\Z[G]}\big(H_1(M_K;\Z[G]),\Z[G])\oplus H_1(M_{C_{\a}};\Z[\Z])\otimes_{\Z[\Z]}\Z[G] ,\Z[G]\big)\\
\cong &\Ext^1_{\Z[G]}\big(H_1(M_K;\Z[G]),\Z[G])\oplus \Ext^1_{\Z[G]}(H_1(M_{C_{\a}};\Z[\Z])\otimes_{\Z[\Z]}\Z[G],\Z[G]\big)\\
\cong&\Ext^1_{\Z[G]}\big(H_1(M_K;\Z[G]),\Z[G])\oplus
\Ext^1_{\Z[\Z]}(H_1(M_{C_{\a}};\Z[\Z]),\Z[\Z]\big)\\
\cong& \Ext^1_{\Z[\Z]}(H_1(M_{C_{\a}};\Z[\Z]),\Z[\Z]\big).\ea \] Note that
$H_1(M_{C_{\a}};\Z[\Z])\cong H_1(S^3\sm C_\a;\Z[\Z])$, in particular it is $\Z$--torsion free.
It follows from \cite[Theorem~3.4]{L77} that
$\Ext^1_{\Z[\Z]}(H_1(M_{C_{\a}};\Z[\Z]),\Z[\Z]\big)\cong
 H_1(M_{C_{\a}};\Z[\Z])$, which is non-trivial since $\Delta_{C_{\a}}(t)\ne 1$. The
same calculation for $\varphi_{\b}$ then gives the contradiction.
\end{proof}

\subsection*{Non-uniqueness of ribbon groups}

A ribbon knot can in general have ribbons  with non-isomorphic groups. For example let $K$ be a non-trivial knot with $\Delta_K(t)=1$, eg, the Whitehead double of a non-trivial knot. Let $G:=\pi_1(S^3\sm K)$, then $G$ is in particular a ribbon group for the knot $L:=K\#-K$.
 On the other hand, $L$ still has trivial Alexander polynomial and so by Theorem~\ref{thm:main} it
follows that $L$ also has a h--ribbon with fundamental group $\Z$. In fact, if $K$ is the Whitehead double of a ribbon knot then it is ribbon, not just h--ribbon, with fundamental group $\Z$ (and so is $L$).


\newpage
\pagenumbers{2501}{2504}   
\count0=2501
{\small\hypertarget{Err}{This} erratum was
received and accepted on 6 September 2006 and published on 18 October
2006 as \href{http://dx.doi.org/10.2140/gt.2006.10.2501}{Geom. Topol. 10 (2006)
2501--2504}.}
\vglue 20pt
\section*{\cl{Correction to `New topologically slice knots'}}
\addcontentsline{toc}{section}{Erratum}
\setcounter{section}0
\setcounter{figure}0
\def\thefigure{\arabic{figure}}
\def\theshorttitle{New topologically slice knots: erratum}

\vglue10pt
 {\leftskip 25pt\rightskip25pt\small 

{\bf Abstract}\qua In \figref{example2} of \cite{FT05} (above) we gave
an incorrect example for Theorem \ref{thm:main}. In this note we
present a correct example.\par} \vglue10pt

 We first recall the satellite construction for knots. Let $K,C$ be
knots. Let $A\subset S^3\setminus K$ be a curve, unknotted in $S^3$.
Then $S^3 \setminus \nu A$ is a solid torus. Now let
$\psi:\partial(\overline{\nu A})\to
\partial(\overline{\nu C})$ be a diffeomorphism which sends a
meridian of $A$ to a longitude of $C$, and a longitude of $A$ to a
meridian of $C$. The space
\[ (S^3\setminus \nu A)\cup_{\psi} (S^3\setminus \nu C) \] is a 3-sphere and the image of $K$ is denoted by
$S=S(K,C,A)$. We say $S$ is the satellite knot with companion $C$,
orbit $K$ and axis $A$. Note that by doing this construction we
replaced a tubular neighborhood of $C$ by a knot in a solid torus,
namely $K\subset S^3\setminus \nu A$.

\begin{figure}[ht!]
\begin{center}
\labellist\small
\pinlabel $K$ at 56 68
\endlabellist
\includegraphics[scale=0.25]{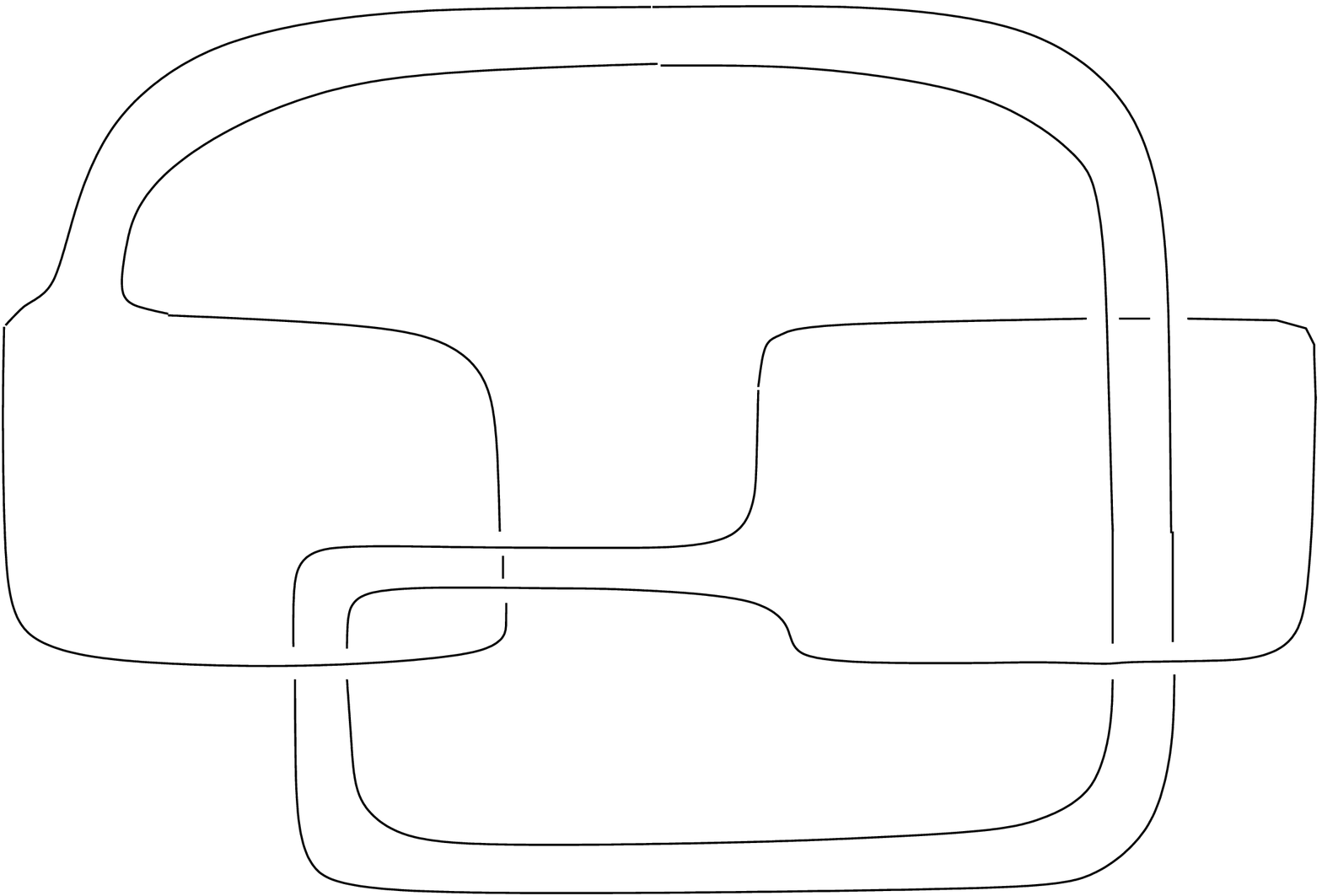}
\caption{The knot $K$} \label{example2noa}
\end{center}
\end{figure}

We now consider the knot $K$ in Figure \ref{example2noa}.
Note that $K$ is the knot $6_1$.  $K$ clearly bounds an
immersed band; pushing this band into $D^4$ we can resolve the
singularities to get a smooth slice disk $D$ for $K$.

It follows from Proposition \ref{propsatribbon} that if $A\subset
S^3\setminus K$ is a curve such that $[A]$ represents an element in
$\ker\{\pi_1(S^3\setminus K)\to \pi_1(D^4\setminus D)\}$, then
$S(K,C,A)$ is topologically slice. We recall that
$\pi_1(D^4\setminus D)$ is isomorphic to the semi--direct product
\[
 \langle a, c \mid aca^{-1}=c^2 \rangle \cong \sr.
\]
Here the generator $a$ of $\Z$ acts on the normal subgroup $\Z[1/2]$
via multiplication by $2$. In \figref{example2} we proposed a
curve $A$ and claimed that it represents the trivial element in $
\pi_1(D^4\setminus D) \cong \sr$. Unfortunately we miscalculated the
image of $A$ in $\sr$. In fact this $A$ represents a non--trivial
element in $\pi_1(D^4\setminus D)$. Hence the curve $A$ of
\figref{example2} does not give an example for
Proposition~\ref{propsatribbon}. We now present a correct example.

Perhaps the first example of a pair $K,A$ which satisfies the above
conditions which comes to mind is to take $K,A$ which form a slice
link $K \cup A$. But it is easy to see that the null--concordance
from $K \cup A$ to a trivial link $K'\cup A'$ induces a concordance
of $S(K,C,A)$ to $S(K',C,A')$. But clearly $S(K',C,A')$ is the trivial
link. This shows that in this case $S(K,C,A)$ is slice. We therefore
have to find examples of $K,A$ such that $K\cup A$ is not slice.

\begin{figure}[ht!]
\begin{center}
\labellist\small\hair 5pt
\pinlabel $K$ at 56 68
\pinlabel $A$ [b] at 281 462
\endlabellist
\includegraphics[scale=0.25]{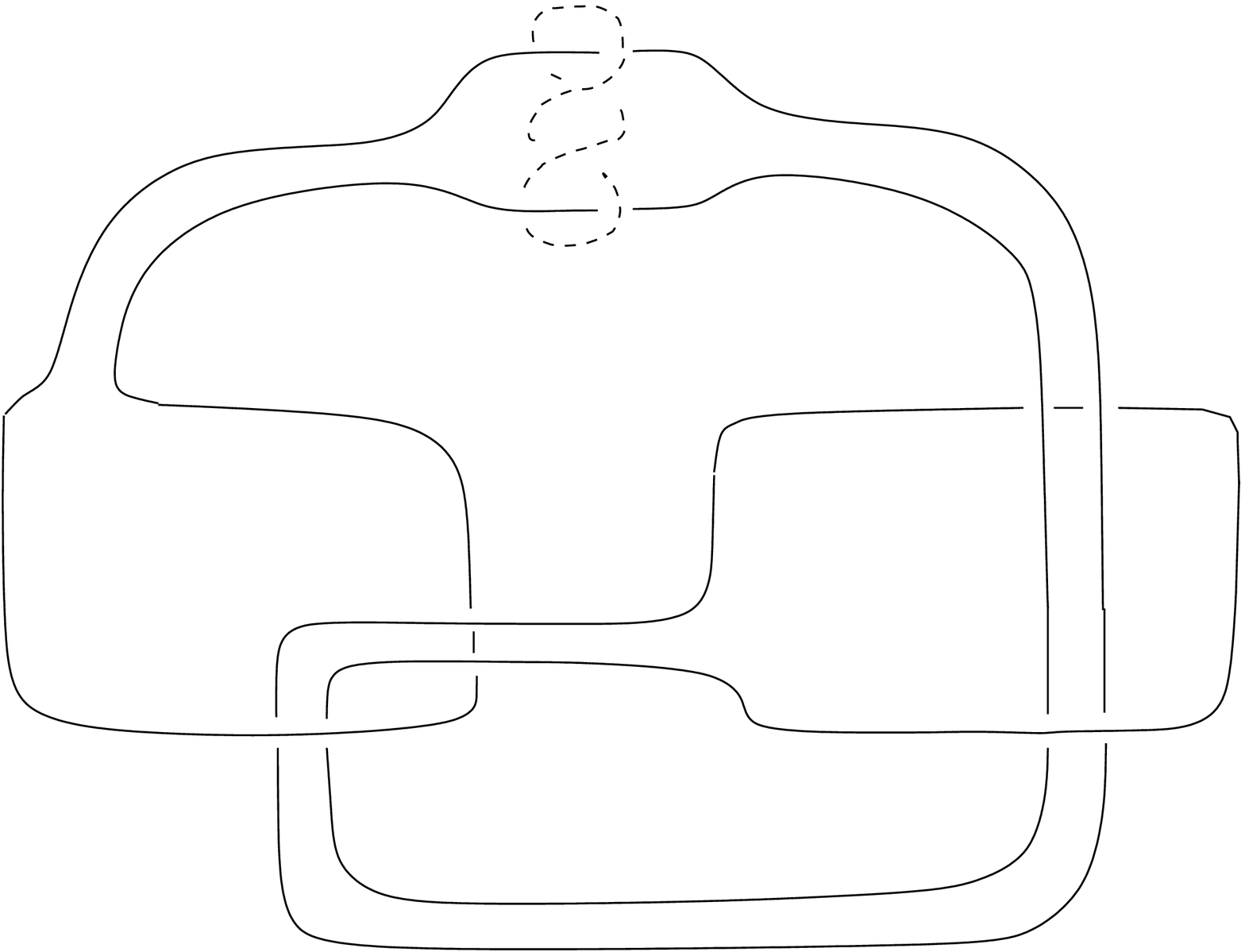}
\caption{The knot $K$ and the curve $A$} \label{example2new}
\end{center}
\end{figure}

Now let $A$ be the  simple closed curve of Figure \ref{example2new}.
Since $D\cap S^3=K$ we can resolve the crossings of $A$ using a
homotopy in $S^3\setminus K\subset D^4\setminus D$. We get a curve
without crossings which is a meridian for the band. Now we push this
curve into $D^4$ `beyond $D$' and then we can contract this curve.
This shows that $A$ is null--homotopic in $D^4\setminus D$. A
straightforward calculation shows that the Alexander polynomial of
the link $K\cup A$ is non--trivial, hence the link $K\cup A$ is not
slice by Kawauchi \cite{Ka78}.

\begin{figure}[h]
\begin{center}
\begin{tabular}{ccccc}
\includegraphics[scale=0.3]{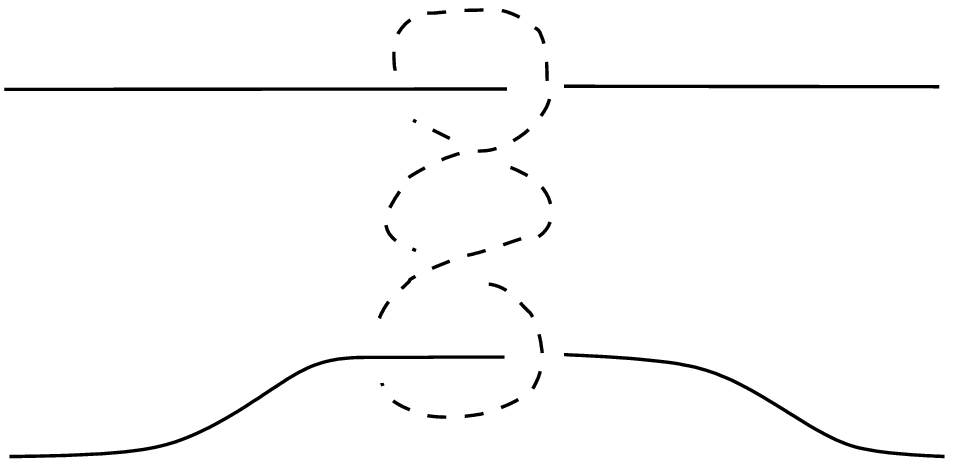}&\hspace{1.5cm}&&&
\includegraphics[scale=0.3]{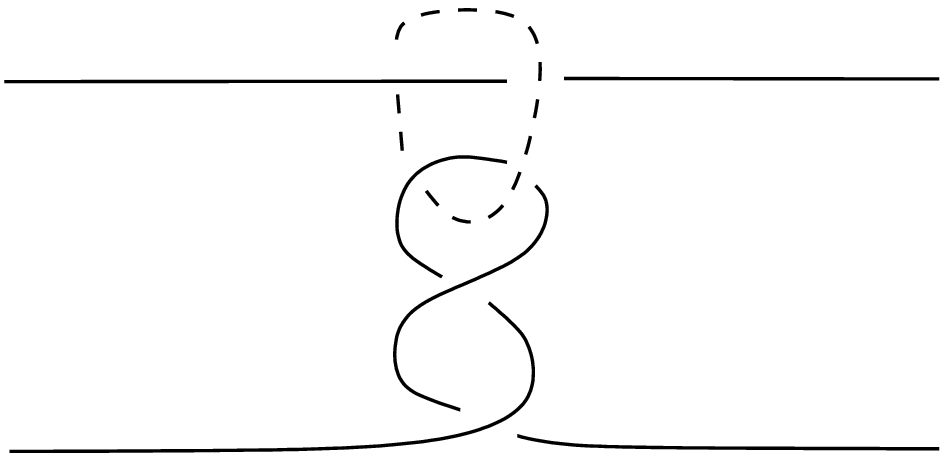}
\end{tabular}
\caption{Untwisting $A$.} \label{fig:twist}
\end{center}
\end{figure}

\begin{figure}[h]
\begin{center}
\begin{tabular}{ccccc}
\includegraphics[scale=0.25]{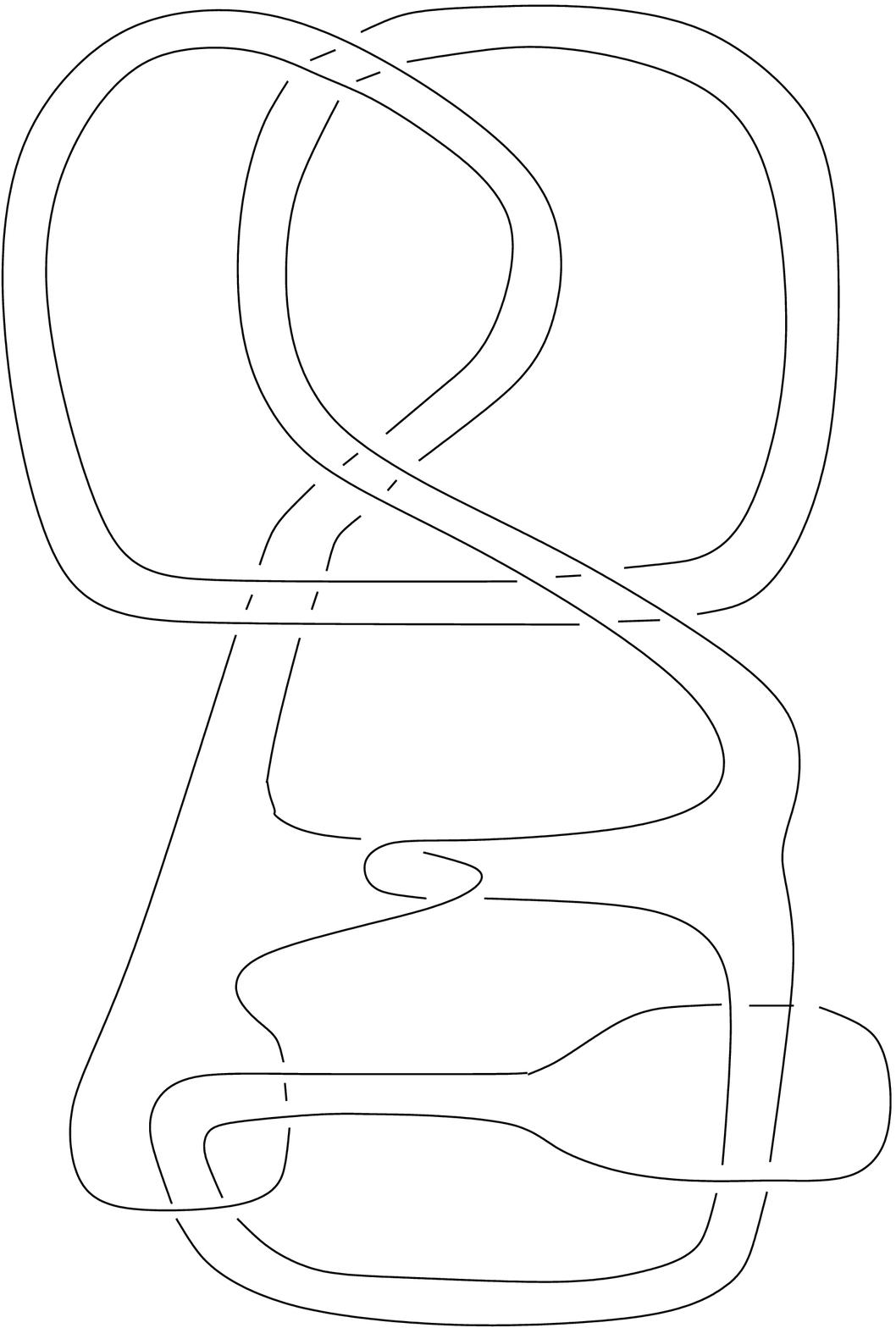}
\end{tabular}
\caption{Satellite knot of the figure-8 knot.} \label{fig:fig8}
\end{center}
\end{figure}

Finally we point out that by untwisting $A$ (and therefore twisting
$K$)  as in Figure \ref{fig:twist} we get a diagram of $K$ in a
`planar' torus. Wrapping this torus around a knot $C$ gives
immediately a diagram for $S(K,C,A)$. For example if we take $C$ to
be the figure-8 knot we get the  diagram in Figure \ref{fig:fig8}
with $26$ crossings.

 We point out that in general if $C$ has a diagram with
crossing number $c$ and writhe $w$, then $S(K,C,A)$ has clearly a
diagram of crossing number $4c+2|w|+10$. This is significantly lower
than the crossing number for the (incorrect) example of $A$ given in
\figref{example2} and will hopefully put our examples within
reach of Rasmussen's $s$--invariant.

\small
{\it Department of Mathematics, Rice University, Houston, TX 77005, USA}\nl
and\nl
{\it Department of Mathematics, University of 
California, Berkeley, CA 94720, USA}

\vglue 15pt
Erratum received: 6 September 2006\nl
Published: 18 October 2006

\end{document}